\documentclass[10pt]{article}
\usepackage{amsmath,amsthm,amscd,amssymb,mathrsfs,setspace}
\usepackage{latexsym,epsf,epsfig}
\usepackage{epstopdf}
\usepackage[hmargin=3cm,vmargin=2cm]{geometry}
\usepackage{color}
\newcommand{\ds}{\displaystyle}

\newcommand{\reals}{\mathbb{R}}
\newcommand{\realstwo}{\mathbb{R}^2}
\newcommand{\realsthree}{\mathbb{R}^3}
\newcommand{\xb}{{\bf{x}}}

\newcommand{\Dx}{{\partial_x}}
\newcommand{\pd}{{\partial}}
\newcommand{\cA}{{\mathscr{A}}}

\newcommand{\Dn}{\partial_{\nu}}

\newcommand{\cE}{{\mathcal{E}}}
\newcommand{\cD}{\mathscr{D}}
\newcommand{\bA}{\mathbb{A}}
\newcommand{\bB}{\mathbb{B}}
\newcommand{\bP}{\mathbb{P}}
\newcommand{\bX}{\mathbb{X}}
\newcommand{\cF}{\mathscr{F}}
\newcommand{\cL}{\mathscr{L}}

\newcommand{\R}{\mathbb{R}}
\newcommand{\cO}{\mathcal{O}}
\newcommand{\cU}{\mathcal{U}}
\newcommand{\Om}{{\Omega}}

\newcommand{\la}{{\lambda}}

\newcommand{\bH}{\mathbf{H}}
\newcommand{\s}{\sigma}
\newcommand{\e}{\epsilon}
\newcommand{\Ez}{E_{z}}
\newcommand{\cH}{\mathcal{H}}

\newcommand{\vr}{{\varrho}}

\newcommand{\tvr}{{\tilde{\varrho}}}
\newcommand{\tv}{{\tilde{v}}}
\newcommand{\tp}{{\tilde{p}}}

\newcommand{\di}{{\rm div\, }}
\newcommand{\g}{{\nabla}}

\theoremstyle{plain}
\newtheorem{theorem}{Theorem}[section]
\newtheorem{lemma}[theorem]{Lemma}
\newtheorem{proposition}[theorem]{Proposition}
\newtheorem{corollary}[theorem]{Corollary}
\newtheorem{assumption}{Assumption}[section]
\theoremstyle{remark}
\newtheorem{remark}{Remark}[section]
\numberwithin{equation}{section}
\numberwithin{theorem}{section}
\numberwithin{remark}{section}

\title{Flow-Plate Interactions:  Well-posedness and Long-Time Behavior}
\date{\today}

 \author{\normalsize \begin{tabular}[t]{c@{\extracolsep{.8em}}c@{\extracolsep{.8em}}c}
           Igor Chueshov  & Irena Lasiecka & Justin T. Webster \\
\it        Kharkov National Univ. & \it Univ. of Virginia    &\it Oregon State University \\
\it        Kharkov, Ukraine & \it Charlottesville, VA &\it Corvallis, OR\\
\it       chueshov@univer.kharkov.ua   &  \it il2v@virginia.edu&  \it websterj@math.oregonstate.edu
\end{tabular}}

\begin{document}
\maketitle

\begin{abstract} {\noindent
We consider flow-structure interactions  modeled  by a modified wave equation coupled at an interface with equations of nonlinear elasticity.  Both subsonic and supersonic flow velocities are treated with Neumann type flow conditions, and a novel treatment of the so called Kutta-Joukowsky flow conditions are given in the subsonic case. The  goal of the paper is threefold: (i) to
  provide an accurate review of recent results
on existence, uniqueness, and stability of weak solutions, (ii) to present a construction of finite dimensional, attracting sets
corresponding to the structural dynamics and discuss convergence of trajectories, and (iii) to state several open questions associated  with the topic.  This second  task is based on a decoupling technique which reduces the analysis of the full flow-structure system to a
PDE system with delay.

\noindent {\bf Key terms}: flow-structure interaction, nonlinear plate,
 nonlinear semigroups, well-posed\-ness, global attractors, PDE with delay. }
\end{abstract}

\section{Introduction}
Flow-structure models have  attracted considerable attention in the past mathematical literature, see, e.g., 
\cite{bal0,bal4,bal3,b-c, LBC96,b-c-1,chuey,springer,jadea12,dowellnon,ACC,webster}
and the references therein.
However, the vast majority of the work done on flow-structure interactions has been devoted to numerical and experimental studies; see,
for instance,
\cite{BA62,bolotin,dowellnon,dowell1,HP02,jfs} and also the survey~\cite{Li03}
and the literature cited there. 
This is not surprising, taking into considerations multitudes of engineering applications with a 
prime example of flutter suppression. 
 Many mathematical studies have been based on linear, two dimensional, plate models  with specific geometries. In these studies  the primary goal was to determine the speed at which flutter occurs
\cite{bal0,BA62,bolotin,dowell1,HP02,Li03}.
 More recently, the study of linear models wherein the two dimensional dynamics are reduced (typical section)
 to a one dimensional structure (beam, or panel with constant width which is infinitely extended in one direction)
with Kutta-Joukowsky flow boundary conditions have enjoyed renewed interest, and
have been extensively pursued in \cite{bal4,shubov3,shubov1}.  This line of work has focused on spectral properties of the system, with particular emphasis on identifying aeroelastic {\em eigenmodes} corresponding to  the associated Possio integral  equation.
\par
In contrast, our interest here concerns PDE aspects of the problem, including the fundamental issue of well-posedness
of {\it  finite energy}  solutions corresponding to nonlinear flow-plate interaction, along with long-time behavior of solutions within the realm of dynamical systems.
\par

In this paper we address  one of the most general flow-structure  PDE models which describes the interactive dynamics between a nonlinear plate and a surrounding potential flow \cite{bolotin,dowellnon}. This class of models is standard in the applied mathematics  literature for the modeling of flow-structure
interactions and dates back to classical discussions  (\cite{bolotin,dowellnon}, and also \cite{B,dowell1} and the references therein).
Owing to the importance  of the associated physical phenomena (including {\em flutter}),
 there has been an immense amount of activity over the last 60 years with the aim of
better understanding  the modeling  and control  aspects  of the underlying dynamics.
An excellent source of recent information on the subject is provided by the review 
paper \cite{B} and also \cite{bal0,bal4,C}.   These include  inspiring  modeling, computational, and experimental studies  of the flow-plate and flow-beam problem.

The main goal of the present paper is to provide PDE based  analysis  of flow-plate dynamics corresponding to large range of  flow velocities  including  both  {\it subsonic and supersonic } flow velocities. As a starting point, we consider  PDE models introduced in
 \cite{dowellnon,dowell1}, and more recently discussed in \cite{B}. A characteristic feature of these models is the absence of the {\it rotational inertia} term. This is in agreement with physics of the problem where the plate is of infinitesimal thickness  and the in-plane accelerations are not accounted for \cite{lagnese}.
 On the other hand, it is by now well known that the presence of rotational terms provides an additional regularizing effect
 on the transverse velocity of the plate, which in turn leads to several desirable mathematical properties such as compactness, gain of derivatives,  etc.. Thus, it was natural that first mathematical PDE  theory of flow-plate dynamics  was carried out under the assumption that rotational inertia is present. With this assumption there is a rather rich theory  that has been developed
 for  both well-posedness and long-time behavior of the  (structural) dynamics subjected to
 {\it mechanical damping} implemented on the plate (see  \cite{springer} and the references therein).

 In view of this, the distinct feature of our work is that the analysis presented herein (a)  {\it does not account for the regularizing effects
 of rotational terms}, and (b) {\it does not {\em require} any mechanical damping imposed on the structure}, 
 which is precisely in agreement with the physical model considered in \cite{dowellnon,B,dowell}.
 Yet, the final discussion presented  in this paper provides   long-time asymptotic  properties of  the originally ``rough'' dynamics
 which ultimately become (without any added damping) smooth and finite dimensional in nature. In order to gain an insight as to
 how it is conceivable to obtain ``attracting behavior" without any dissipation present, we offer some preliminary explanation:
experimentally we see that the flow (particularly at supersonic speeds) has the ability of inducing stability in the moving structure. This is the case when the structure itself does not possess any mechanical or frictional damping. This dissipative effect is not immediately noticeable in the standard energy balance equation. However, a reduction technique introduced in \cite{LBC96,b-c-1} allows us to write the full flow-structure interaction as a certain delayed plate model, and demonstrates the stabilizing effects of the flow, provided that {\em rotational terms are not present in the model}. The full flow dynamics manifest themselves in the form of non-conservative forces acting upon the structure via the {\em downwash} of the flow. In the case where rotational inertia is present in the plate model, the downwash of the flow is not substantial enough to dissipate the mass term due to inertia.
Thus, the challenge to resolve this paradox and to solve the mathematical  problem of long-time behavior for the model is to show that the non-conservative  effect of the downwash leads to the desired dissipative long-time dynamics. The obstacles are then the lack of compactness and the criticality of nonlinear term (again, due to the absence of rotational inertia), along with the lack of gradient structure in the reduced (delayed) dynamical system. 

As compactness of the dynamics and gradient structure are pillars of long-time behavior analysis,
previous studies have only been successful when regularizing (compactifying) effects are present in the model.
 Hence, the long-time behavior of subsonic and supersonic models without rotational 
 inertia---or other regularizing terms (i.e. thermoelastic plates) in the dynamics \cite{ryz,ryz2}---has 
 been virtually unchartered territory. The tools required to tackle the mathematical issues of loss of compactness, 
 loss of gradient structure in dynamical systems, and supercriticality of nonlinear terms in the 
 equations have been developed only recently.  
 This paper presents the state of the art results and  mathematical methods as are relevant to the problem at hand.
 For other aspects of the problem  such as computational and experimental we refer reader to a 
 large wealth of literature available at present (see the discussion and the references in
 \cite{bal0,BA62,bolotin,A,Li03}).

\subsection{Notation}
For the remainder of the text we write $\xb$ for $(x,y,z) \in \realsthree_+$ or $(x,y) \in \Omega \subset \realstwo_{\{(x,y)\}}$, as dictated by context. Norms $||\cdot||$ are taken to be $L_2(D)$ for the domain $D$. The symbols $\nu$ and $\tau$ will be used to denote the unit normal and unit tangent vectors do a given domain, again, dictated by context. Inner products in $L_2(\realsthree_+)$ are written $(\cdot,\cdot)$, while inner products in $L_2(\R^2\equiv\pd\R^3_+)$ are written $<\cdot,\cdot>$. Also, $ H^s(D)$ will denote the Sobolev space of order $s$, defined on a domain $D$, and $H^s_0(D)$ denotes the closure of $C_0^{\infty}(D)$ in the $H^s(D)$ norm
which we denote by $\|\cdot\|_{H^s(D)}$ or $\|\cdot\|_{s,D}$. We make use of the standard notation for the trace of functions defined on $\realsthree_+$, i.e. for $\phi \in H^1(\realsthree_+)$, $\gamma[\phi]=\phi \big|_{z=0}$ is the trace of $\phi$ on the plane $\{\xb:z=0\}$. (We use analogous notation for $\gamma[w]$ as the trace map from $\Omega$ to $\partial \Omega$.)
\subsection{Basic Principles of the Fluid/Gas Dynamics}
The model in consideration describes the interaction between a nonlinear plate and a fluid or flow of gas above it.
Suppose  the domain $\cO= \R^3_+$ is filled with fluid whose dynamics
are governed by the {\em compressible} Navier--Stokes
 (or Poisson--Stokes) system (see, e.g., \cite{chorin-marsden}, \cite{dowell1} and \cite{landau}
for the physical backgrounds of the models)
 which are written for the density $\tvr$,  velocity $\tv$ and
pressure $\tp$:
\begin{align*}
 &
   \tvr_t+\di\{\tvr \tv\}=0\quad {\rm in\quad} \cO
   \times(0,+\infty), 
   \\[2mm]
 &
   \tvr\big[\tv_t+ (\tv,\g) \tv)\big] =\di T(\tv,\tp) \quad {\rm in\quad} \cO
   \times(0,+\infty), 
\end{align*}
In this section we utilize the notation $\xb=(x_1;x_2;x_3)$ to accommodate the tensor analysis, and
\[
\left[\di T(v,p)\right]^i=\sum_{j=1}^3\pd_{x_j} T^{ij}(v,p), \quad i=1,2,3.
\]
where  $T=\{T^{ij}\}_{i,j=1}^3$ is the stress tensor of the fluid,
\[
T^{ij}\equiv T^{ij}(v,p)=\mu_f\left(v^i_{x_j}+v^j_{x_i}\right)+\big[\la_f\, \di v-p\big]\delta_{ij}.
\quad i,j=1,2,3.
\]
Here $\mu_f$ and $\la_f$ are (non-negative) viscosity coefficients (which  vanish
in the case of  invisid fluid).
We assume that the fluid is isothermal, i.e., the pressure $\tp$ is a \emph{linear} function of the density
$\tvr$. For simplicity we take $\tp=\tvr$.
\par
Now we consider  linearization of the model
with respect to some reference state $\{\vr_*;v_*;p_*\}$. We suppose
that unperturbed flow  $v_*$
represents the fluid/gas moving with speed $U$ along the $x_1$-direction, i.e.,
$v_*=(U,0,0)$  and $\vr_*,p_*$ are constants.
For simplicity we assume $p_*=\vr_*=1$. Then
for small perturbations $\{\vr; v; p\}$ we obtain the equations:
\begin{subequations}\label{bar-model1-U}
\begin{align}
 &
   \vr_t+ U\vr_{x_1}+\di\, v=0\quad {\rm in}~~ \cO
   \times\R_+,\label{den-eq1U}
   \\[2mm]
 &
   v_t+ U v_{x_1} -\mu_f\Delta v -(\mu_f+\la_f)\g \di\, v + \g\vr =0 \quad {\rm in}~ \cO
   \times\R_+,\label{flu-eq1U}
\end{align}
\end{subequations}
We need to also supply \eqref{bar-model1-U} with appropriate boundary conditions.
We choose impermeable boundary conditions, which look as follows:
\begin{equation*}
(T\nu,\tau)=0~~{\rm on}~\pd\cO,~~~ (v,\nu)=0~~{\rm on}~\pd\cO\setminus\Om,~~ (v,\nu)=u_t+Uu_{x}~~{\rm on}~\Om,
\end{equation*}
where $\nu$ is the unit outer normal to $\partial\cO$, $\tau$ is a unit tangent direction to $\partial\cO$, $u$ is the deflection
of flexible part $\Om$ of the boundary (see some discussion in \cite{chorin-marsden} and
also in \cite{dowell1} where an explanation of the term $u_t+Uu_{x}$
on the boundary is given).
In the case  $\cO= \R^3_+$, the boundary conditions have the form
\begin{equation}\label{impermeability2}
\mu_f\left(v^i_{x_3}+v^3_{x_i}\right)=0,
~~ i=1,2,~~{\rm on}~\pd\cO,~~~ v^3=0~~{\rm on}~\pd\cO\setminus\Om,~~ v^3= u_t+Uu_{x}~~{\rm on}~\Om,
\end{equation}

This model  describes the case
of (possibly viscous) {\em compressible} gas/fluid flows
and was recently studied in \cite{Chu2013-comp} in the viscous
case $\mu_f>0$ with the {\em zero} speed
$U$ of unperturbed flow. This general model also
leads to several special cases which are important from an applied
point of view.
\begin{itemize}

 \item \textbf{Incompressible}  fluid, i.e. $\di \, v=0$ and $\vr$ constant:
 In the
    \emph{viscous} case $\mu_f>0$ the  standard linearized Navier-Stokes equations arise; fluid plate interaction
   in this case were studied in  \cite{ChuRyz2011,cr-full-karman,ChuRyz2012-pois,berlin11}. Results
    on well-posedness and attractors
   for different plate situations and domains were obtained.
 The \emph{invisid} case was studied in \cite{Chu2013-invisid}  in the same context.
  \item
\textbf{Invisid  compressible} fluid ($\mu_f=0$ and  $\la_f=0$): In this case from \eqref{bar-model1-U}
we can obtain wave type dynamics    for the (perturbed) velocity potential $(v=-\g \phi$, potential flow) of the form
(see also   \cite{BA62,bolotin,dowell1}):
\begin{equation}\label{flow}\begin{cases}
(\partial_t+U\partial_x)^2\phi=\Delta \phi & \text { in } \realsthree_+ \times (0,T),\\
\phi(0)=\phi_0;~~\phi_t(0)=\phi_1,\\
\Dn \phi = -\big[(\partial_t+U\partial_x)u (\xb)\big]\cdot \mathbf{1}_{\Omega}(\xb)& \text{ on } \realstwo_{\{(x,y)\}} \times (0,T).
\end{cases}
\end{equation}
In these variables the pressure/density of the fluid has the form $\vr=(\partial_t+U\partial_x)\phi$.
Due to \eqref{impermeability2} in the case of the perfect fluid ($\mu_f=0$ and  $\la_f=0$)
we have only the one boundary condition given above.
This is exactly the model for a gas/fluid on which we concentrate in this paper.
\end{itemize}

\subsection{PDE Description of the Full Gas-Plate Model}
 The behavior of the plate is governed by the second order (in time) Kirchoff plate equation; 
 we consider the von Karman nonlinearity, which is the principal 
 `physical' nonlinearity used in the modeling of the large oscillations of thin, 
 flexible plates---so called \textit{large deflection theory} \cite{ciarlet,lagnese}.
\par
The environment we consider is $\realsthree_+=\{(x,y,z): z \ge 0\}$. The plate
 is modeled by  a bounded domain $\Omega \subset \reals^2_{\{(x,y)\}}=\{(x,y,z): z = 0\}$ with smooth boundary $\partial \Omega = \Gamma$.
 The plate is  embedded  (perhaps partially) in  a large rigid body (producing 
 the  so called \textit{clamped} boundary conditions) immersed in an inviscid 
 flow (over body) with velocity $U \neq 1$ in the  negative $x$-direction.
  We have normalized $U=1$ to be Mach 1, i.e. $0 \le U <1$ corresponding to subsonic flows and $U>1$ is supersonic.
\par
The scalar function $u: \Om\times \R_+ \to \reals$ represents the vertical displacement of the plate in the $z$-direction at the point $(x;y)$ at the moment $t$. We take the plate with clamped or free-clamped boundary conditions:

\begin{equation}\label{plate0}\begin{cases}
u_{tt}+\Delta^2u+f(u)= p(\xb,t) & \text { in } \Om\times (0,T),\\
 u(0)=u_0;~~u_t(0)=u_1,
\end{cases}
\end{equation}
Clamped boundary conditions are given by:
\begin{equation}\label{cl}
u=\Dn u = 0 , \text{ on } \pd\Om\times (0,T).
\end{equation}
Free-clamped boundary conditions are given by:
\begin{align}\label{free}
u=\Dn u = 0 ~ \text{ on } ~ \pd\Om_1\times (0,T),\nonumber \\
\mathcal{ B}_1 u =0,   \mathcal{B}_2 u =0 ~ \text{on}~ \pd\Om_2\times (0,T).
\end{align}
where $\pd\Om_i,~ i =1,2$ are complementary parts of the boundary $\partial \Omega $ , and $\mathcal{B}_1$, $ \mathcal{B}_2$ represent moments and shear forces, given by \cite{lagnese}:
\begin{eqnarray*} 
\mathcal{B}_1 u \equiv \Delta u + (1-\mu) B_1 u ~
\notag \\
\mathcal{B}_2 u \equiv \partial_{\nu} \Delta u + (1-\mu) B_2 u - \mu_1 u -
\beta u^3;  ~~  \beta \geq 0.
\end{eqnarray*}
The  boundary operators $B_1$ and $B_2$ are given by \cite{lagnese}:
\begin{equation*}
\begin{array}{c}
B_1u = 2 \nu_1\nu_2 u_{xy} - \nu_1^2 u_{yy} - \nu_2^2 u_{xx}\;, \\
\\
B_2u = \partial_{\tau} \left[ \left( \nu_1^2 -\nu_2^2\right) u_{xy} + \nu_1
\nu_2 \left( u_{yy} - u_{xx}\right)\right],
\end{array}
\end{equation*}
where here $\nu=(\nu_1, \nu_2)$ is the outer normal to $\Gamma$, $\tau= (-\nu_2,
\nu_1)$ is the oriented unit tangent vector along $\Gamma$. (Note that in the case of $\realsthree_+$ with boundary $\realstwo$, $\nu$ is identified with the $-z$ direction.) The parameters $\mu_1$
and $\beta $ are nonnegative, with the constant $0<\mu<1$ having the meaning of the
Poisson modulus.

The choice of boundary conditions is dependent upon the application. For instance, panels mounted in a rigid frame will satisfy clamped boundary conditions. Alternatively, certain flag-like configurations in air flow provide a classical example of free-clamped combination. (See Section \ref{config} below.)

The aerodynamical  pressure  $p(\xb,t)$  represents by the third component $T^{33}$
of  the stress tensor describes the coupling with the flow
and will be given below.
In this paper we consider the typical nonlinear (cubic-type) force term resulting from  aeroelasticity modeling \cite{bolotin,dowell,dowell1,HolMar78}: the {\em von Karman model}, given by $$f(u)=-[u, v(u)+F_0],$$ where $F_0$ is a given function
  from $H^{3+\delta}(\Om)$. In the case of free-clamped boundary conditions we must also assume that $F\in H^2_0(\Om)$. 
    
  The von Karman bracket $[u,v]$  is given by
\begin{equation*}
[u,v] = \partial ^{2}_{x} u\cdot \partial ^{2}_y v +
\partial ^{2}_y u\cdot \partial ^{2}_{x} v -
2\cdot \partial_x\partial_y u\cdot \partial_x\partial_y
v,
\end{equation*} and
the Airy stress function $v(u,w) $ solves the following  elliptic
problem
\begin{equation*}
\Delta^2 v(u,w)+[u,w] =0 ~~{\rm in}~~  \Omega,\quad \Dn v(u,w) = v(u,w) =0 ~~{\rm on}~~  \pd\Om
\end{equation*}
 (we make use of the notation $v(u)=v(u,u)$).
Von  Karman equations are well known in nonlinear elasticity and
constitute a basic model describing nonlinear oscillations of a
plate accounting for  large displacements, see \cite{karman}
and also \cite{springer,ciarlet}  and
references therein.
\par

By taking the aerodynamical pressure
of the form
\begin{equation}\label{aero-dyn-pr}
p(\xb,t)=p_0+\big(\partial_t+U\partial_x\big)\gamma[\phi]
\end{equation}
 in (\ref{plate0}) above, $p_0 \in L_2(\Omega)$,
 and  using  the perturbed wave equation in \eqref{flow}, we arrive at
  the fully coupled model:
\begin{equation}\label{flowplate}\begin{cases}
u_{tt}+\Delta^2u+f(u)= p_0+\big(\partial_t+U\partial_x\big)\gamma[\phi] & \text { in } \Om\times (0,T),\\
u(0)=u_0;~~u_t(0)=u_1,\\
BC(u)& \text{ on } \pd\Om\times (0,T),\\
(\partial_t+U\partial_x)^2\phi=\Delta \phi & \text { in } \realsthree_+ \times (0,T),\\
\phi(0)=\phi_0;~~\phi_t(0)=\phi_1,\\
\Dn \phi = -\big[(\partial_t+U\partial_x)u (\xb)\big]\cdot \mathbf{1}_{\Omega}(\xb) & \text{ on } \realstwo_{\{(x,y)\}} \times (0,T).
\end{cases}
\end{equation}
where $BC(u)$ represent appropriate plate boundary conditions as specified in  (\ref{free}).
\begin{remark}
The flow-structure interaction models which we discuss here, established in the literature,  do not involve moving frames (boundaries) of reference. The model corresponds to a ``snap shot"  of the structure interacting with the flow at a given time $t$ with $\Omega$ as a fixed frame of reference in the plane.  This is in contrast with analyses where the interface between the fluid and solid is intrinsically evolving in time \cite{ref6,aiim}.
\end{remark}

\subsection{Physical Configurations and Other Flow Boundary Conditions}\label{config}
At this point, we pause to point out other modeling configurations of great interest. The  mathematical arguments  presented below are done in the context of purely {\em clamped} boundary conditions, i.e. no portion of the boundary $\partial \Omega$ is free. We have presented the model above in additional generality, including the possibility of a free component, however, the full treatment of a free portion of the boundary is technically challenging and currently under investigation.  This is due to the  loss of sufficient  regularity of the boundary data imposed for the flow. Clamped boundary conditions
assumed on the boundary of the plate allow for smooth extensions to $\R^2$ of the Neumann flow boundary conditions
satisfied by the flow. In the absence of these, one needs to approximate the original dynamics in order to construct
 sufficiently smooth functions amenable to PDE calculations.  This is a technical challenge and was carried out in a similar fashion in \cite{supersonic}, although the need for this analysis was not due to plate boundary conditions.  Since free boundary conditions   are physically interesting, and represent challenging well-posedness problems,  we discuss other possible configurations of the model and their physical pertinence in this section.
\begin{remark}
So called {\em hinged} or {\em pinned} (homogeneous) boundary conditions do not differ substantially (mathematically) from the clamped boundary conditions. (See \cite{springer,lagnese} for details.) These boundary conditions do not yield a markedly different analysis than that of clamped conditions, with the key point being that in both cases the plate component of the model admits a well-posed biharmonic problem, and hence the Sobolev $H^2(\Omega)$ norm is equivalent to the norm $||\Delta \cdot||$. It is of interest to note that boundary control via moments are typically expressed through a hinged type condition. See \cite{springer} for details in the case of the plate alone, and \cite{webster&lasiecka} for a well-posedness analysis of these boundary conditions in the flow-plate model.
\end{remark}
Perhaps the most mathematically interesting (and difficult) case corresponding to this model is the 
{\em free-clamped} configuration. These types of boundary conditions are extremely important in the 
modeling of airfoils and in the modeling of panels in which some component of the boundary is left free.
In addition to the difficulties associated with the free plate boundary condition coupled 
with a Neumann type flow condition on the interface (the last line of \eqref{flow}), another 
natural {\em flow} boundary condition is the so called {\em Kutta-Joukowsky} (K-J) condition 
(a dynamic, mixed-type boundary condition) \cite{bal0,K1,dowell,K2,K3} and considered below. 
The applicability of this flow boundary condition is highly dependent upon the geometry of the plate in question. Physically, the K-J conditions corresponds to taking a zero pressure jump in the flow {\em off of the wing}---specifically at the trailing edge \cite{dowell1}. Mathematically, we are extremely interested in understanding the affect of 
the K-J condition interacting with both clamped and free type plate boundary conditions.

The so called {\em Kutta-Joukowsky} conditions specify that (for the flow): \begin{equation}\label{KJ} \begin{cases}
\partial_{\nu} \phi = -(\partial_t+U\partial_x)u &\text{on} ~\Omega\\
\psi=\phi_t+U\phi_x = 0 &\text{off} ~\Omega.\end{cases}
\end{equation}
These mixed type flow boundary conditions are taken to be accurate for plates in the clamped-free configuration. 
A recent analysis was made of these boundary conditions coupled with a clamped plate; 
we include an outline of this analysis in Section \ref{KJsec}. 
However, more in depth analysis of the arguments presented reveals that 
free boundary conditions can be incorporated as well.

The configuration below represents an attempt to model oscillations of a plate which is {\em mostly free}. The dynamic nature of the flow conditions correspond to the fact that the interaction of the plate and flow is no longer static along the free edge(s), and in this case the implementation of the K-J condition is called for. In this case we take the free-clamped plate boundary conditions, and the mixed flow boundary conditions:
\begin{equation}\label{physical}\begin{cases}
u=\Dn u = 0,  &\text{ on }  \pd\Om_1\times (0,T) \\
\mathcal{B}_1 u =0,   \mathcal{B}_2 u =0,~ &\text{on}~ \pd\Om_2\times (0,T)\\
\Dn \phi = - (\partial_t+U\partial_x)u, ~& \text{on}~\Omega\times (0,T)\\
\Dn \phi= 0, & ~\text{on}~\Theta_1\times(0,T)\\
\psi=\phi_t+U\phi_x=0, &~ \text{on}~\Theta_2\times(0,T)\end{cases}
\end{equation}
with the regions described by the image below (the regions $\Theta_i$ extend in the natural way into 
the remainder of the $x-y$ plane):

\begin{center}
\includegraphics[scale=.23]{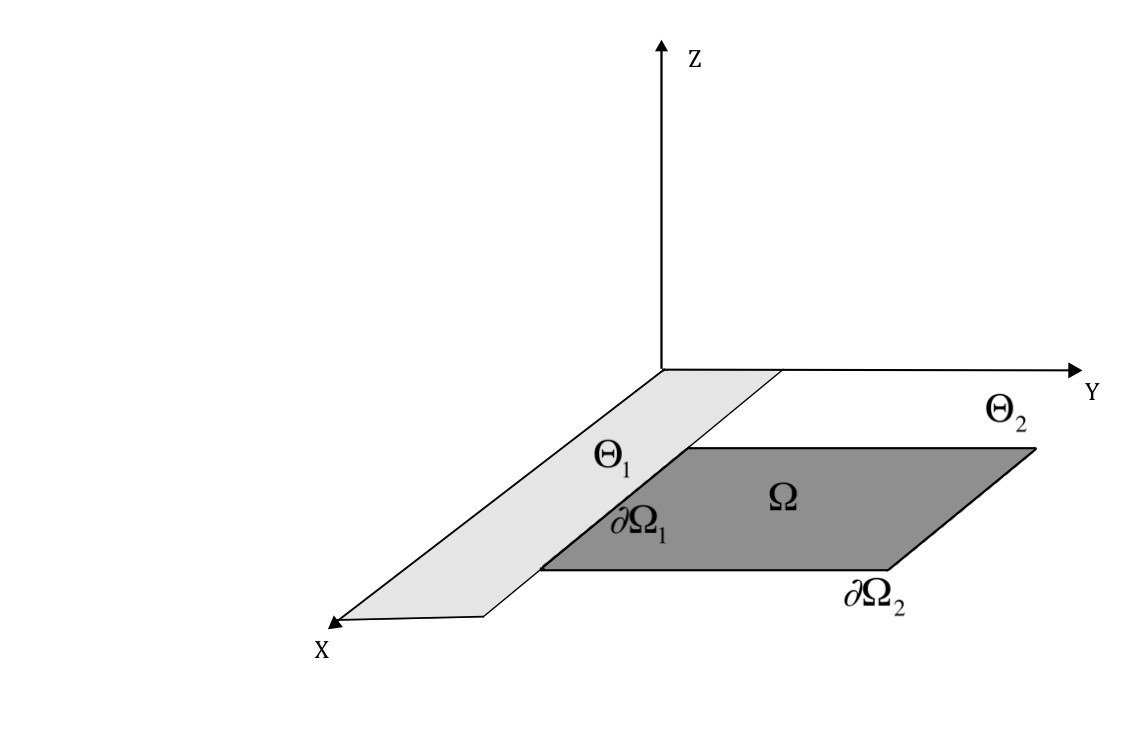} \end{center}

\begin{remark}
The configuration above arises in the study of airfoils, but another related configuration referred to as {\em axial flow} takes the flow to occur in the $y$ direction in the picture above. In our analysis, the geometry of the plate (and hence the orientation of the flow) do not play a central role. In practice, the orientation can have a dramatic effect on the occurrence and magnitude of the oscillations associated with the flow-structure coupling. In the case of axial flow, the above configuration is often discussed in the context of {\em flag flutter} or flapping. See \cite{K4} for more details. 
\end{remark}

The physical nature of the model given above makes its analysis desirable; however, this model involves a high degree of mathematical complexity due to the dynamic and mixed nature of the boundary coupling near the flow-plate interface. From the point of view of the existing analysis, much of the well-posedness and long-time behavior analysis 
benefits upon the ability to extend the plate solution by 0 outside of $\Omega$.
This is particularly true with Neumann bounary conditions imposed on the flow. 
Such extension is  not a trivial task for {\em any} boundary condition other than that of clamped, and hence, considering the analysis herein with mixed boundary conditions is immediately non-trivial. There is some indication in preliminary studies that the K-J flow condition may be more amenable to supporting free-clamped plate boundary conditions. This seems to agree with the implementation of the K-J boundary conditions in other recent analyses \cite{bal0,bal4} (we provide more discussion of this in the Open Problems section, \# 1).
\smallskip\par
\noindent{\em Assumption:}  The analysis to follow is done strictly in the context of the {\em clamped plate}. This is done to streamline exposition and avoid involved technicalities; in what follows, we do specify what issues arise when considering other plate boundary conditions. The analysis for the clamped plate represents a first step---complete with a full spectrum of challenging PDE issues---in producing well-posedness and long-time behavior results for this class of fluid-structure models.

\section{Energies and Phase Spaces}
\label{energiessection}

There is a  distinct difference in the form of energy functional for the {\it subsonic} case and {\it supersonic} case taken with the standard Neumann type flow conditions.
Indeed, the standard energy function for the flow loses ellipticity in the supersonic case. For this reason these two cases are treated separately. We will also consider the subsonic model with the K-J flow conditions below, but its functional setup is similar to the supersonic case, so we discuss its properties following that of the supersonic model.

\subsection{Subsonic Energy}\label{subenergy}
We begin with the {\em simpler} subsonic case. In this case ($0 \le U <1$) energies 
can be derived by applying standard plate and flow multipliers, $u_t$ and $\phi_t$ respectively, 
along with the boundary conditions to obtain the energy relations for the plate and the flow. 
This procedure leads to the energy which is bounded from below. (See \cite{springer,jadea12,webster} for details.) 
We discuss this briefly now.

The description of the system in the subsonic case is given by ($0 \le U<1$):
\begin{equation}\label{syssub}\begin{cases}

(\partial_t+U\partial_x)^2\phi =\Delta \phi& \text{ in } \realsthree_+ \times(0,T),
\\
\Dn \phi = -\big[(\partial_t+U\partial_x)u (\xb)\big]\cdot \mathbf{1}_{\Omega}(\xb) & \text{ on }
 \realstwo_{\{(x,y)\}} \times (0,T),
\\  u_{tt}+\Delta^2u+f(u)=p_0+\gamma[\psi]&\text{ in } \Om\times (0,T),\\
BC(u)& \text{ on } \pd\Om\times (0,T),\\
\end{cases}
\end{equation} where $BC(u)$ correspond to clamped or free-clamped boundary conditions as given above.

In this case, we have  {\it the flow energy} and {\it interactive } energies given by
\begin{align*}
\hat E_{fl}(t) = & \dfrac{1}{2}\big(||\phi_t||^2-U^2||\partial_x\phi||^2+||\nabla \phi||^2\big),~~
\hat E_{int}(t)=U<\gamma[\phi],u_x>
\end{align*}
The  {\it plate energy} is defined as usual 
\begin{align}\label{energies}
E_{pl}(t) =& \dfrac{1}{2}\big(||u_t||^2+ a(u,u) \big) +\Pi(u),
\end{align}
where $ a(u,u) $ is a bilinear form topologically equivalent to $ H^2(\Omega) $ defined 
by\footnote{In the case of clamped or hinged boundary conditions the bilinear form $a(u,w) $ collapses  \cite{springer} to the familiar expression
$\ds a(u,w) = \int_{\Omega} \Delta u \Delta w dx $.}

\begin{equation*}
a(u,w)=\int_{\Omega }\widetilde{a}(u,w)+\mu _{1}\int_{\Gamma _{1}}uw,
\end{equation*}
with
\begin{equation*}
\widetilde{a}(u,w)\equiv u_{xx}w_{xx}+u_{yy}w_{yy}+\mu
(u_{xx}w_{yy}+u_{yy}w_{xx})+2(1-\mu )u_{xy}w_{xy}.  \label{a-tild}
\end{equation*}

\noindent $\Pi(u)$ is a potential of the nonlinear forces  given
by\footnote{In the case of clamped boundary conditions (i.e. $\partial \Omega_2=\emptyset$) the boundary terms vanish.}
 \begin{equation*}
 \Pi(u)=\frac{1}{4}||\Delta v(u)||^2-\frac12 <[u,u],F_0>+\int_{\pd\Om_2}\left[ \frac{\mu_1}2u^2+\frac{\beta}2 u^4\right] 
 \end{equation*}
 in the case of the von Karman nonlinearity.

The total energy  in the subsonic case is defined as a sum of the three components $$\hat {\cE}(t) =\hat E_{fl}(t) + E_{pl}(t) + \hat E_{int}(t)  $$ and satisfies (formally)
 $\hat{\cE}(t) = \hat{\cE}(s)$.

 \subsection{Supersonic Energy}
 We note, that in the {\it supersonic}   case  ($ U > 1 $)  the flow part of the energy $\hat{E}_{fl}(t) $ is no longer nonnegative.
 This, being  the source of major mathematical difficulties, necessitates  a  different approach.
In the case of {\em supersonic} flows we make use of the flow acceleration multiplier $$(\partial_t+U\partial_x)\phi \equiv \psi.$$ We thus have a change of ``state" from the natural (hyperbolic type) state variable associated to the flow equation: $\phi_t \rightarrow (\phi_t+U\phi_x) = \psi$.
\par
We then have a description of our  coupled system as follows:
\begin{equation}\label{sys1}\begin{cases}

(\partial_t+U\partial_x)\phi =\psi& \text{ in } \realsthree_+ \times(0,T),
\\
(\partial_t+U\partial_x)\psi=\Delta \phi&\text{ in } \realsthree_+\times (0,T),
\\
\Dn \phi = -\big[(\partial_t+U\partial_x)u (\xb)\big]\cdot \mathbf{1}_{\Omega}(\xb) & \text{ on }
 \realstwo_{\{(x,y)\}} \times (0,T),
\\  u_{tt}+\Delta^2u+f(u)=p_0+\gamma[\psi]&\text{ in } \Om\times (0,T),\\
BC(u)& \text{ on } \pd\Om\times (0,T),\\
\end{cases}
\end{equation} where $BC(u)$ correspond to clamped, or free, or free-clamped boundary conditions as given above,
see \eqref{cl} and \eqref{free}.

This leads to the following (formal) energies, arrived at via Green's Theorem:
\begin{align}\label{energies-n}
E_{pl}(t) =& \dfrac{1}{2}\big(||u_t||^2+ a(u,u) \big) +\Pi(u),~~
E_{fl}(t) =  \dfrac{1}{2}\big(||\psi||^2+||\nabla \phi||^2\big),
~~
\cE(t) = E_{pl}(t)+E_{fl}(t),\notag
\end{align}
where $ a(u,u) $ and $ \Pi(u) $ were introduced before.
\par
With these energies, we have the formal energy relation
\begin{equation}\label{energyrelation}
\cE(t)+ U\int_0^t <u_x,\gamma[\psi]> dt = \cE(0). 
\end{equation}
This energy relation reveals an apparent {\it loss of dissipation} and the associated {\it loss of gradient structure}.
It also  provides the first motivation for viewing the dynamics (under our change of phase variable) as comprised of a generating piece and a ``perturbation".
\par
Finite energy constraints manifest themselves in the natural requirements on the functions $\phi$ and $u$:
\begin{equation*}
\phi(\xb,t) \in C(0,T; H^1(\realsthree_+))\cap C^1(0,T;L_2(\realsthree_+)), 
\end{equation*} 
\begin{equation*}
u(\xb,t) \in C(0,T; \cH)\cap C^1(0,T;L_2(\Omega)), 
\end{equation*}
where 

$$
\cH=H_{\pd\Om_1}^2(\Omega)\equiv\left\{ u\in H^2(\Om)\, :\, u=\Dn u = 0 , \text{\rm on } \pd\Om_1  \right\}\footnote{We make use of this notation throughout the remainder of the text where pertinent in order to keep the statements as general as possible. However, for most of the remainder of the paper we specify to the case of clamped boundary conditions.}
$$

\par
In working with well-posedness considerations (and thus dynamical systems), the above finite energy constraints lead to the  so called finite energy space, which we will take as our state space: 
\begin{equation}\label{space-Y}
Y = Y_{fl} \times Y_{pl} \equiv \big(H^1(\realsthree_+)\times L_2(\realsthree_+)\big) \times \big(\cH \times L_2(\Omega)\big).
\end{equation}

 The new representation of the energies as in (\ref{energies}) provides good topological measure for the sought after solution, however the energy balance is {\it lost} in (\ref{energyrelation}) and, in addition, the boundary term is
 ``leaking energy" and  involves the traces of $ L_2$ solutions, which are possibly {\it not defined}  at all.
 
 \subsection{Subsonic Case Taken with the K-J Flow Condition}
In this section we briefly relate the functional setup above (supersonic with Neumann flow conditions) with the case of subsonic flows taken with the K-J flow condition given in \eqref{KJ}. Our full system in this case is given by
\begin{equation}\label{flowplate2}\begin{cases}
u_{tt}+\Delta^2u+f(u)= p_0+\gamma[\psi] & \text { in } \Om\times (0,T),\\
u(0)=u_0;~~u_t(0)=u_1,\\
u=\Dn u=0 & \text{ on } \pd\Om\times (0,T),\\
(\partial_t+U\partial_x)\phi=\psi & \text { in } \realsthree_+ \times (0,T),\\
(\partial_t+U\partial_x)\psi=\Delta \phi & \text { in } \realsthree_+ \times (0,T),\\
\phi(0)=\phi_0;~~\phi_t(0)=\phi_1,\\
\Dn \phi = -(\partial_t+U\partial_x)u & \text{ on } \Omega   \times (0,T).\\
(\partial_t+U\partial_x)  \phi = 0 & \text{ on } \realstwo \backslash \Omega \times (0,T).\end{cases}
\end{equation}

In line with the analysis of the supersonic flows (with standard flow boundary conditions) in \cite{supersonic}, we make use of the flow acceleration multiplier $(\partial_t+U\partial_x)\phi \equiv \psi$ for the {\em subsonic flow}, taken with the K-J flow boundary conditions. 
Thus for the flow dynamics, instead of $(\phi;\phi_t)$ we again have
the state variables  $(\phi;\psi)$.
This leads to the same (formal) energies, arrived at via Green's Theorem when the multipliers $u_t$ and $\psi$ are applied to \eqref{flowplate} as in the supersonic case described above.
These energies provide the formal energy relation for the system (implementing the K-J conditions)
\begin{equation}\label{energyrelationkj}
\cE(t)+ U\int_0^t <u_x,\gamma[\psi]> dt = \cE(0). \end{equation}

As in the supersonic case, the relations above lead to the same finite energy space.

\section{Main Results}
The main results presented below address: (i) the well-posedness of finite energy solutions, 
(ii) the construction of an attracting set for the dynamics of the structure, 
and (iii) convergence of the full flow-plate dynamics to stationary states for a certain class of configurations. 
All results presented below  implement {\em clamped boundary conditions for the plate.}

\subsection{Definitions of Solutions}\label{solutions}
In the analysis we will encounter strong (classical), generalized (mild), and weak (variational) solutions.
In the well-posedness analysis found in \cite{supersonic} semigroup theory is used, which requires 
the use of \textit{generalized} solutions; these are strong limits of strong solutions.
These solutions satisfy an integral formulation of (\ref{flowplate}), and are also called \textit{mild} by some authors.
We note that generalized solutions are also weak solutions, see, e.g., \cite[Section 6.5.5]{springer} and \cite{webster}.
Weak solutions possess the quality of {\em finite energy} and satisfy the corresponding variational relation.
We do not  provide here rigorous definitions of these solutions and instead refer to \cite{supersonic} and \cite{webster};
see also \cite{springer} for detailed discussion of solutions for different types of boundary conditions in the case of stand-alone plate models.
All of the solutions mentioned below are generalized, and differ only their smoothness properties.
In all statements below we try to indicate clearly the corresponding smoothness. 

\subsection{Well-Posedness of the Model}
We now present an overview of the well-posedness results for the flow-plate model.
	\subsubsection{Statement of Subsonic Results}
	
	 \label{sec:WellP:Subsonic:Clamped0:Results}
We assume that $0\le U <1$ (subsonic) and $p_0 \in L_2(\Omega)$, with $F_0 \in H^{3+\delta}(\Omega)$.

\begin{theorem}[Linear System]\label{linearresult} Suppose $f \equiv 0$.
\begin{enumerate}
\item \textbf{Linear Generation} The dynamics operator associated to \eqref{syssub} generates a strongly continuous semigroup $T_t$  on the state space $$Y=H^1(\realsthree_+) \times L_2(\realsthree_+) \times H_0^2(\Omega)\times L_2(\Omega).$$
\item \textbf{Strong Solutions} Assume
$$
u_1 \in H^2_0(\Omega), ~u_0\in W \equiv\{w\in H_0^2(\Omega) : \Delta ^2 w \in L_2(\Omega)\}= H^4(\Omega)\cap H^2_0(\Omega).
$$ 
Moroever, suppose $$\phi_0 \in H^2(\realsthree_+)~\text{ and }~\phi_1 \in H^1(\realsthree_+)$$ 
with the following compatibility condition in place: 
\begin{equation*}
\partial_z \phi_0\big|_{z=0} = \begin{cases}u_1+U\partial_x u_0 & \text{ if } \xb \in \Omega\\ 0 & \text{ if } 
\xb \notin \Omega \end{cases}.
\end{equation*} 
Then \eqref{syssub} has a unique strong solution for any interval $[0,T]$. 
This solution possesses the properties 
\begin{equation*} (\phi;\phi_t;\phi_{tt}) \in L_{\infty}(0,T;H^2(\realsthree_+)\times H^1(\realsthree_+)\times L_2(\realsthree_+)),
\end{equation*}
\begin{equation*}
(u;u_t;u_{tt})\in L_{\infty}(0,T;W \times H_0^2(\Omega) \times L_2(\Omega))
\end{equation*} 
and satisfies the energy equality $\hat{\mathcal{E}}(t)=\hat{\mathcal{E}}(s)$ for $t>s$, with $\hat{\mathcal{E}}(t)$ defined as in Section \ref{energiessection}.
\item \textbf{Generalized and Weak Solutions} Assume 
$$u_0 \in H_0^2(\Omega),~ u_1 \in L_2(\Omega),~ \phi_0 \in H^1(\realsthree_+),~ 
\text{ and } ~\phi_1 \in L_2(\realsthree_+).$$ 
Then \eqref{syssub} has a unique generalized solution on any interval $[0,T]$. 
This solution satisfies $\hat{\mathcal{E}}(t) = \hat{\mathcal{E}}(s)$ for $t>s$ which then implies
that there exist constants $ C >0, \omega > 0 $ such that
\begin{equation*}
||T_t||_{\cL(Y) } \leq C e^{\omega t }, t > 0
\end{equation*}
 Every generalized solution is also weak.
\end{enumerate}
\end{theorem}

\begin{theorem}[{\bf Nonlinear Semigroup}]
\label{vonkarmansolution} If $f(u)=-[v(u)+F_0,u]$, then under the \textbf{strong solutions} assumptions Theorem
\ref{linearresult}, for all $~T$, \eqref{syssub} has a unique strong solution on $[0,T]$; under the \textbf{generalized and weak solutions} assumptions, for all  $~T$, \eqref{syssub} has a unique generalized (and hence weak) solution on $[0,T]$ denoted by $S_t (y_0) $. This is to say that $(S_t, Y) $ is a (nonlinear) dynamical
system on $Y$. Moreover, this solution is uniformly bounded in time in the norm of the state space $Y$.
This means that there exists constant C such that $$ ||S_t (y_0)||_Y \leq C (||y_0||_Y).$$
\end{theorem}

\begin{remark}
We note that while the solution to linear problem does not need to  be bounded in time,
the corresponding solution to the nonlinear problem is ultimately bounded.
This is due to the presence of the nonlinearity which controls the interactive part of the energy  at the level of higher frequencies. \end{remark}
\begin{remark} This theorem remains true if (a) the von Karman nonlinearity is replaced by any nonlinearity $f: H^2(\Omega) \cap H_0^1(\Omega) \to L_2(\Omega)$ which is locally Lipschitz (e.g. the Berger or Kirchoff type nonlinearities \cite{supersonic}), and (b) the nonlinear trajectories possess a uniform bound in time.
\end{remark}

	\subsubsection{Statement of Supersonic Results}\label{supres}
	\begin{theorem}[{\bf Linear}]\label{th:lin}
Consider linear problem in (\ref{flowplate}) with $f(u) =0$.  Let $ T > 0 $. Then,   for every initial datum
$( \phi_0, \phi_1; u_0, u_1 ) \in Y $ there exists unique generalized solution
\begin{equation}\label{phi-u0reg}
(\phi (t), \phi_t(t); u(t), u_t(t)) \in C([0, T ],  Y)
\end{equation}
 which   defines  a $C_0$-semigroup
$T_t : Y \rightarrow Y $ associated with  (\ref{sys1}) (where $ f =0$).
\par
 For any initial data in
  \begin{equation}\label{Y1}
Y_1 \equiv \left\{ y= (\phi,\phi_1;u,u_1)  \in Y\;   \left| \begin{array}{l}
\phi_1 \in H^1(\R^3_{+} ),~~ u_1 \in \cH,  ~\\ -U^2 \Dx^2 \phi  + \Delta \phi \in L_2(\R^3_{+}), \\
\Dn  \phi = - [u_1 +U \Dx u ]\cdot {\bf 1}_{\Omega} \in H^1(\R^2), \\
~
-\Delta^2 u + U\gamma[ \Dx \phi] \in L_2(\Omega) \end{array} \right. \right\}
\end{equation}
   the corresponding  solution is also strong
in the sense that 
\begin{equation*}
(\phi (t), \phi_t(t); u(t), u_t(t)) \in C([0, T ],  Y_1)
\end{equation*}
   \end{theorem}
 We shall turn next to nonlinear problem.
\begin{theorem}[{\bf Nonlinear Semigroup}]\label{th:nonlin}
Let $ T > 0 $ and let $f(u)$ be given by  the von Karman nonlinearity. Then, for every initial datum
$( \phi_0, \phi_1; u_0, u_1 ) \in Y $ there exists unique generalized solution
$(\phi, \phi_t; u, u_t)$ to (\ref{flowplate})
with the clamped boundary conditions
 possessing property \eqref{phi-u0reg}. This solution is also weak and generates
  a  nonlinear  continuous  semigroup
$S_t : Y \rightarrow Y$ associated with  (\ref{sys1}).
\par
If $( \phi_0, \phi_1; u_0, u_1 ) \in Y_1$, where
 $Y_1 \subset Y$ is given by \eqref{Y1}, then the corresponding solution is also strong.
\end{theorem}
 \begin{remark}
  In comparing the results obtained with a subsonic case, there are two major differences at the qualitative level:

First, the regularity of strong solutions obtained in the subsonic case \cite{jadea12,webster} coincides with regularity expected for classical solutions. In the supersonic case, there is a  loss of differentiability in the flow  in the tangential  $x$ direction, which then propagates to the  loss of differentiability in the structural variable $u$. As a consequence strong solutions do not exhibit sufficient regularity in order to perform the  needed calculations.
  To cope with the problem, special regularization procedure was introduced in \cite{supersonic}
  where  strong solutions are approximated by sufficiently regular functions though not solutions to the given PDE.
  The limit passage allows to obtain the needed estimates valid for the original solutions \cite{supersonic}.

 Secondly, in the subsonic case one shows that the resulting semigroup is {\it bounded} in time,
  see \cite[Proposition 6.5.7]{springer} and also \cite{jadea12,webster}.
  This property  could not be shown in the supersonic analysis, and most likely does not hold.  The leak of energy in the energy relation can not be compensated for by the nonlinear terms (unlike in the subsonic case).
   \end{remark}

\subsubsection{Discussion of  Well-posedness Results for the Subsonic Case taken with the Kutta-Joukowsky Flow Conditions}
Some aspects of both the subsonic and supersonic analyses discussed in the previous theorems appear in the analysis of the case of subsonic flows with the K-J flow condition. However, new technical ingredients enter in a substantial way. 
They depend critically on  the {\it subsonic } range of velocities for the flow (unlike the Neumann boundary conditions). 
The work on well-posedness of  finite energy solutions  for this system is very recent, and as such, 
we include an outline of the proof in Section \ref{KJsec}. 

We suffice to say, at this point, that in Section \ref{KJsec} we demonstrate well-posedness of mild solutions to the system given in \eqref{flowplate2} on the finite energy space $Y$; this result depends upon a critical assumption (Assumption \ref{hilb2} in Section \ref{tracereg}) which then yields additional trace regularity of the acceleration potential $(\partial_t+U\partial_x)\phi \big|_{\realstwo}$. Such an assumption is not needed for the well-posedness analysis in the case of the standard Neumann type data (Theorems \ref{vonkarmansolution} and \ref{th:nonlin} above). Owing to the rather technical nature of this assumption, which deals with the invertibility of a two dimensional singular integral transform (related to the so called Possio equation), we do not state the result here. Rather, we provide a self contained discussion and outline of the well-posedness proof in Section \ref{KJ}. Additionally, we demonstrate that this key assumption is satisfied  when the dimensionality of the model is reduced to a two dimensional structure interacting with a one dimensional flow. 
Due to the reduction of dimension, the said Assumption \ref{hilb2} is reduced to invertibility of {\it finite} Hilbert transforms on $L_p$, $ p \in (1,2) $. This fact was first observed and used in \cite{bal0} where linear spectral analysis is used in order to determine  flutter modes.

\subsection{Long-time Behavior and Attracting Sets}
We present here the  long-time characteristics of the dynamical system induced by ($S_t, Y$)
in both subsonic and supersonic cases with the standard Neumann flow condition.
As noted in the introduction, the problem at hand is rendered challenging for the following reasons:
\begin{itemize} \item
 the {\it unboundedness of the flow domain};
 \item the {\it  lack of inherent dissipation}
for both flow and plate;
\item the {\it  lack of compactness} in both the flow and plate dynamics;
\item the {\it absence of gradient structure} for the dynamical system.
\end{itemize}

Nevertheless,  the result we obtain demonstrates that the trajectories for the structure are attracted by a finite dimensional and smooth attracting set.  In order to obtain this type of result, a certain amount of decoupling of the flow dynamics from structure dynamics  is absolutely necessary. This will be accomplished by realizing that the flow effects produce (after some time) {\em delayed forces}
 acting upon the plate. Thus, instead of considering the full dynamical system with the nice properties of invariance, we will need to consider a modified dynamical system for which desirable properties can be shown {\em for the structure}.

  \begin{theorem}\label{th:main2}
Suppose $0\le U \ne 1$,  $F_0 \in H^{3+\delta}(\Omega)$ and $p_0 \in L_2(\Omega)$.
 Then there exists a compact set $\mathscr{U} \subset \cH \times L_2(\Omega)$ of finite fractal dimension such that $$\lim_{t\to\infty} d_{\cH \times L_2(\Omega)} \big( (u(t),u_t(t)),\mathscr U\big)=\lim_{t \to \infty}\inf_{(\nu_0,\nu_1) \in \mathscr U} \big( ||u(t)-\nu_0||_2^2+||u_t(t)-\nu_1||^2\big)=0$$
for any weak solution $(u,u_t;\phi,\phi_t)$ to (\ref{flowplate}) taken with clamped boundary conditions ($\cH=H_0^2(\Omega)$), and with
initial data
$$
(u_0, u_1;\phi_0,\phi_1) \in  \cH\times L_2(\Omega)\times H^1(\realsthree_+)\times L_2(\realsthree_+)
$$
which are
localized  in $\R_+^3$ (i.e., $\phi_0(\xb)=\phi_1(\xb)=0$ for $|\xb|>R$ for some $R>0$). Additionally we have the extra regularity $\mathscr{U} \subset \big(H^4(\Omega)\cap \cH \big) \times H^2(\Omega)$.
\end{theorem}

While Theorem \ref{th:main2} pertains to ``attractors" capturing the limiting dynamics of the plate component of solutions, 
it is of interest to develop
long-time behavior results describing the entire evolution, if possible. 
In particular, it is of interest (and often important) to determine when the end behavior of solutions is simple, 
perhaps {\em static} (in contrast to potentially chaotic behavior). 
The physical phenomenon associated with static end behavior for the plate is known as {\em buckling}
(large deflection theory \cite{lagnese}). 
The issue which we address in the theorem below is under which circumstances  full trajectories (plate and flow) 
converge to stationary states---solutions to the static problem corresponding to \eqref{syssub}.
We may implement  additional damping   of the 
form $ku_t$ with $k>0$ in order to strengthen the result of Theorem \ref{th:main2} in this direction.

\begin{theorem}\label{conequil}
Let $0\le U<1$ and suppose the mechanical damping term $ku_t$, $k>0$ is incorporated into \eqref{flowplate}
above on the LHS of the equation; assume $p_0 \in L_2(\Omega)$ and $F_0 \in H^{3+\delta}(\Omega)$. Then any weak solution $(u(t),\phi(t))$ to the system with localized initial flow data (as in Theorem \ref{th:main2}) has the property that $$\lim_{t \to +\infty}\inf_{\bar u, \bar \phi \in \mathcal N} \left\{||u(t)-\bar u||^2_{H^2(\Omega)}+||u_t(t)||^2_{L_2(\Omega)}+||\phi(t)-\bar \phi||_{H^1(B_{\rho})}^2+||\phi_t(t)||^2_{L_2(B_{\rho})} \right\}=0,$$ where $\mathcal N$ is the set of stationary solutions to \eqref{flowplate}. 
\end{theorem}

The proof of Theorem \ref{conequil} which makes critical use of (i) mechanical damping imposed on the structure, 
(ii) compactness of the attractor stated above, 
and (iii) the finiteness of the dissipation integral
$$
k \int_0^{\infty}||u_t||_{L_2(\Omega)}^2 d\tau < \infty
$$ 
via the energy relation in Section \ref{subenergy}. For some details see \cite{supersonic}
and also \cite[Section 12.4.2]{springer}.

\begin{remark}
As recently documented in \cite{jfs} the result of Theorem \ref{conequil} is not valid in the supersonic case, even with a small supersonic speeds. It  turns out that any amount of  structural damping is not a position to eliminate flutter in 
low modes for low supersonic speeds. In that sense the result claimed in Theorem \ref{conequil} is optimal. 
\end{remark}

The proof of Theorem \ref{th:main2} will be discussed later in Section 6. 

\subsection{Discussion of Results in Relation to Past Literature}

\subsubsection{Well-posedness}
For the reasons described above, well-posedness results in the past literature dealt mainly with the dynamics possessing some regularizing effects. This has been accomplished by either accounting for non-negligible rotational forces \cite{b-c,b-c-1,springer}  or  strong damping from thermal effects \cite{ryz,ryz2}.  In the first case, the  linear part of plate dynamics becomes
\begin{equation}\label{alpha}
u_{tt} - \alpha \Delta u_{tt} + \Delta^2 u = p_0+(\phi_t + U \phi_x)\big|_{\Om};
\end{equation}
while in the second case, thermal effects lead to the consideration of a strongly damped plate
\begin{equation}\label{alpha1}
u_{tt}  -\alpha \Delta u_{t} + \Delta^2 u =  p_0+(\phi_t + U \phi_x)\big|_{\Om}.
\end{equation}
In both cases the plate velocity has the property $u_t \in H^1(\Omega) $---which is the needed regularity
for the applicability (in this scenario) of many standard tools in nonlinear analysis. Even at this stage, though, the problem is far from simple. One is faced with low regularity of boundary traces,
due to the failure of Lopatinski conditions (unlike the Dirichlet case for the flow, where there is no loss of regularity to  wave  solutions in their boundary traces).
In fact, the first contribution to the mathematical analysis of the problem  is \cite{LBC96,b-c-1} (see also \cite[Section 6.6]{springer}), where the case  $\alpha > 0$ (rotational) is fully treated.  The method employed in \cite{LBC96,b-c-1,springer}
relies on the following main ingredients:
(i) sharp microlocal estimates for the solution to the wave equation driven by $H^{1/2}(\Omega) $ Neumann  boundary data given by $u_t + U u_x$. This gives $\phi_t|_{\Omega} \in
L_2(0,T; H^{-1/2}(\Omega))$  \cite{miyatake1973} (in fact more regularity is presently known:
$H^{-1/3}(\Omega) $ \cite{tataru});  and (ii) the regularizing effects on the velocity $u_t$ (i.e. $u_t \in H^1(\Omega)$) due to the fact that  $\alpha >0$.
 The above items, along with an explicit solver for the three dimensional perturbed wave equation, and a Galerkin approximation for the structure, allow one to construct  a fixed point for the appropriate
 solution map. The method  works well in both cases $0<U<1$ and  $U > 1$.  Similar ideas were used  more recently \cite{ryz,ryz2} in  the case when thermoelastic effects are added to the model; in this case the dynamics also exhibit
$H^1(\Omega)$ regularity of the velocity in both the rotational and non-rotational cases due to the analytic regularizing  effects induced by thermoelasticity \cite{redbook}.

However, when $\alpha  =0$ (non-rotational case), and additional smoothing is not present, $u_t$ is only $L_2(\Omega)$.  In that case the corresponding estimates
 become singular, destroying the  applicability of the  previous methods.  Naturally,  the first analysis of the problem {\em with} $\alpha=0$ (no rotational inertia) depends critically on the condition $U <1$. The  main  mathematical difficulty of this problem is  the presence of  the   boundary  terms:  $ (\phi_t +
 U \phi_x)|_{\Omega} $  acting as the aerodynamical pressure on the plate. When $U =0$, the corresponding  boundary terms exhibit monotone behavior with respect to   the energy  inner product
 (see \cite[Section 6.2]{springer} and \cite{cbms})
 which is topologically equivalent  to  the topology of the  space $Y$ given by \eqref{space-Y}.   The latter  enables  the use of monotone operator theory (\cite[Section 6.2]{springer} and \cite{cbms}).  However, when
  $U > 0$ this is no longer true  with respect to the  topology induced by the energy spaces.  The  lack of the natural dissipativity for both
interface traces, as well as the nonlinear term in the  plate equation,
make the  task of proving well-posedness challenging.

A recent book \cite[Chapter 6]{springer} provides an account of relevant results in the case $\alpha\geq 0$ and $0 \le U <1$. Existence of a nonlinear semigroup capturing finite energy solutions   has been shown for the first time  in \cite{webster}, and later expanded upon in  \cite{jadea12}; 
the proof of Theorem \ref{linearresult} given in \cite{webster} relies on two  principal ingredients: 
(i) renormalization of the state space which yields $\omega$-dissipativity 
for the dynamics operator (which is nondissipative in the standard norm on the state space); 
(ii) the sharp regularity of Airy's stress function, see \cite{fhlt} and also  \cite[Corollary 1.4.5, p. 44]{springer}, which converts a supercritical nonlinearity into a critical one
\begin{equation*}
||v(u)||_{W^{2,\infty}(\Omega) } \leq C ||u||^2_{H^2(\Omega)};
\end{equation*}
and (iii) control of low frequencies for the system by the nonlinear source, represented by the inequality
\begin{equation}\label{potentialbound}
||u||_{2-\delta}^2\leq \epsilon [ ||u||^2_{H^2(\Omega)} + ||\Delta v(u)||^2_{L_2(\Omega) } ] + C_{\epsilon},~~\forall ~\epsilon>0, ~~\delta \in (0,2] .
\end{equation}
The above inequality holds for clamped  and simply supported  boundary conditions.
In the case of
 clamped-free conditions (\ref{free}) we need to assume that $\beta > 0$ and to add the boundary energy term. 
\par
The  issue of well-posedness in the presence of {\it supersonic flows} was, to a certain extent, the final open case for this class of standard models with Neumann flow data, and was not effectively handled by previous theories. It has only been fully resolved recently in \cite{supersonic}, with the main results  reported above in Section \ref{supres}.
 The method and results in \cite{supersonic}
{\it  do not depend on any smoothing mechanism} (as we take $\alpha =0$), and they are applicable {\em for all} $U \ne 1$.
The key ingredients rely on the development of a suitable {\it  microlocal  trace theory}  for the velocity of the flow, and  the implementation of the corresponding  estimates with semigroup theory in extended spaces.
In this way, a-priori estimates allow for a construction of a nonlinear semiflow which evolves finite energy solutions -with corresponding results stated in Theorem \ref{th:nonlin}.

To the knowledge of the authors the models introduced in \eqref{flowplate2} and \eqref{physical} 
have not been addressed in the case of two dimensional plates and three dimensional flows. However, 
in his pioneering work, the author of \cite{bal0,bal4} (and references therein) considers a 
linear {\em airfoil}  immersed in a subsonic flow; the wing is taken to have 
a high aspect ratio thereby allowing for the suppression of the span variable, 
and reducing the analysis to individual chords normal to the span. By reducing the problem to a one dimensional analysis, many technical hangups are avoided, and Fourier-Laplace analysis is greatly  effective.
Ultimately, the problem of well-posedness and $L_p$ regularity of solutions can be realized in the context of the classical {\em Possio integral problem} \cite{bal0,bal4,tricomi}, involving the inverse Hilbert transform and analysis of Mikhlin multipliers. In the approach presented here, we attempt to characterize our solution by similar means and 
point out how the two 
dimensional analysis greatly complicates matters and gives rise to singular integrals in higher dimensions. 
In addition, our analysis pertains to solvability of the full system 
(flow and plate)-rather than solvability of aeroelastic Possio's equation -a critical component, 
however characterizing only boundary interaction of the flow-plate variables.

\subsubsection{Long-time Behavior}
While long-time behavior of von Karman evolutions  has been well studied \cite{springer}, the long-time behavior of both subsonic and supersonic models {\it without}  rotational inertia or other regularizing terms has been virtually uncharted territory. (For the results corresponding to $\alpha > 0 $, see \cite{springer,ryz,ryz2} and references therein.)  One of the obvious reasons is that the tools needed to tackle mathematical issues such as loss of compactness, loss of gradient structure in dynamical system, and supercriticality of nonlinear terms in the equations have been developed only recently.

The first key idea in this analysis is based on reduction of the full flow-plate dynamics to a delayed plate system where flow information is encapsulated in a delay term.
This approach was already  applied to Berger plate models in \cite{oldchueshov1,oldchueshov2}
 for the proof  of  the existence of attractors for the associated reduced plate system with a delay term.
In fact, an abstract  long-time behavior analysis of second order nonlinear PDEs with delays has been developed in
\cite{Chu92b} (see also \cite{springer}): first, in the case of the von Karman model \textit{with} rotational inertia,
and secondly, in  \cite{oldchueshov1,oldchueshov2}, in the case of the Berger model with a ``small" delay term (corresponding to a hypersonic speeds $U>>1$).
 These expositions address the existence and properties of global attractors for this general plate with delay in the presence of a natural form of interior damping and additional regularizations.

We again emphasize that the presence of rotational inertia parameter $\alpha>0$ in the plate model, while drastically improving the topological properties of solutions, is neither natural nor physically desirable in the context of flow-structure interaction.
 Instead, when the rotational inertia term is neglected ($\alpha=0$), the damping secured by the \textit{flow alone} (via the reduction step) provides the main mechanism for stabilization. Additionally, the principal nonlinearity in the analysis of clamped and clamped-free plates (panels) is that of von Karman; the Berger nonlinearity is taken to be an {\em approximation} of the von Karman model based on simplifying physical assumptions \cite{berger0}. 
 In view of this, it is paramount to the problem at hand to consider the long-time behavior of the von Karman plate model {\it  in the absence of rotational inertia} and impose {\it no limiting regimes}  on the  flow velocity parameter $U \geq 0$, $U \ne 1$.

   The mathematical difficulties which arise in this model force us to consider new long-time behavior technologies applied within this framework. Specifically, the abstract approach mentioned above (and utilized in previous long-time analyses for rotational inertia) {\em does not apply} in this case.  We make use of a relatively new technique \cite{chlJDE04,ch-l,springer} which allows us to address the asymptotic compactness property for the associated dynamical system without making reference to any gradient structure of the dynamics (not available in this model, owing to the dispersive flow term). In addition, we are able to demonstrate extra regularity of the attractor via a quasistability approach  \cite{springer}. In proving finite dimensionality and smoothness of the attractors, the criticality of the nonlinearity and the lack of gradient structure prevents one from using
a powerful technique of {\it backward smoothness of trajectories}  \cite{Babin-Vishik,springer,kh}, where smoothness is propagated forward from the equilibria.  Without a gradient structure, the attractor may have complicated structure (not being characterized by the equilibria points). In order to cope with this issue, we take the advantage
of novel method that is based on suitable  approximations  and exploits only the compactness of the attractor.

In the absence of the damping on the flow, it is reasonable to expect that the entire evolution is strongly stable in the sense of strong convergence to the  equilibria.
In fact, such results have been shown for {\it subsonic velocities and plate models 
with rotational forces and additional mechanical damping}  \cite{springer} 
or plates with {\it thermal effects accounted for} \cite{ryz,ryz2}. 
The key to obtaining such results is a {\em viable energy relation} and a priori bounds on solutions which yield {\em finiteness of the dissipation integral}. Since we have seen that the energy relation in the supersonic case has non-dissipative terms which potentially ``leak" energy, we focus on the case of subsonic flows. Additionally, since we must make use of the energy relation coming from the full flow-plate system we may not utilize the damping which was previously ``culled" from the flow. This indicates that some structural damping is necessary to obtain Theorem \ref{conequil}.
 On the other hand recent results reported in \cite{jfs} indicate that such results 
 are false for a low supersonic speeds.
\subsection{Open Problems}
In what follows below we shall list some open problems which naturally emerge from the material
presented above.
\begin{enumerate}
\item
As pointed out in Section 1.3, the  {\it free-clamped } boundary conditions imposed on the plate are of great physical interest.
From the mathematical point of view, the difficulty arises at the level of linear theory  when one attempts to construct ``smooth" solutions of the corresponding evolution. The typical procedure of extending  plate solutions by zero
outside $\Omega$ leads to a jump in the  Neumann  boundary conditions  imposed on the flow.
In order to contend with this issue, regularization  procedures are needed in order obtain smooth approximations
of the original solutions. While  some regularizations have been already introduced  in \cite{supersonic},
more study is needed in order to demonstrate the effectiveness  of this ``smoothing" for the  large array of problems described in this work.  Free boundary conditions in the context of piston theory and boundary dissipation have been recently studied in \cite{bociu}.
\item
We note that the arguments used in the case of K-J boundary conditions below are perhaps more amenable to the consideration of free boundary conditions imposed on the portion of the boundary of the plate. This is in line with certain engineering applications (e.g.  flag type models \cite{K1}). 
On the other hand, existence of attracting sets   in the above configuration is a totally uncharted and very challenging area. In addition,  rigorous analysis of well-posedness for the case of {\em supersonic flows} ($U>1$), taken with the K-J boundary conditions, represents an entirely wide-open problem for any plate boundary conditions. 

\item
While Theorem \ref{th:main2} asserts an existence of compact attracting set with an additional regularity of {\it one derivative}, it would be interesting to see whether for regular loads $p$ and $F_0$ one obtaind $C^{\infty}$
regularity of the attracting set. This will involve additional boot-strap arguments in line with the strategy outlined in \cite{springer} and utilized recently in \cite{pelin}.
\item Recent techniques have shown that geometrically constrained (localized) interior damping is a viable form of dissipation for the von Karman and Berger evolutions \cite{pelin}. Additionally, a preliminary study of the flow-plate interaction (taken with subsonic flows and standard Neumann type coupling) has shown that nonlinear boundary damping is not a viable candidate for stabilizability to a global attractor \cite{ACC}. As such, it is a natural question to ask whether localized interior damping will be a successful control mechanism to establish the existence of a global attractor in the case of a {\em weak coupling} between the flow and the plate, and perhaps produce the analogous result to Theorem \ref{conequil} making use of damping which is active only on a ``small" portion of the domain. 
\item
The treatment presented above excludes the  so called {\em transonic barrier} $ U =1$.
Indeed, for $U=1$ the analysis provided breaks down in the essential way, namely the principal spatial operator associated to the flow $\Delta-U^2\partial_x^2$ becomes degenerate in the $x$ direction. In the case of supersonic or subsonic flows, we are able to exploit the definite sign of this operator associated to the $x$ variable, however, when $U=1$ the flow equations becomes of degenerate hyperbolic type. In addition, according to \cite{dowell,dowell1}, near $U=1$ the flow equation acquires the term $\phi_x\phi_{xx}$, which is referred to the {\em local Mach number} effect. Thus near the transonic barrier the flow equation becomes degenerate and quasilinear. While numerical/experimental work  predicts appearance of shock waves \cite{C}, to our best knowledge  no mathematical treatment of this problem is  available at present.
 \item
 Finally, the ultimate goal is to consider a fully   nonlinear flow model. Experimental-numerical results predicting shock waves  in the evolution  are partially
 available \cite{C}. However the mathematical aspects of   this problem are  presently wide open.

\end{enumerate}
\section{Proof of Well-Posedness: Subsonic Flows with Kutta-Joukowsky Conditions}\label{KJsec}
We see above that the primary source of mathematical difficulty lies in interpreting and unravelling
the energy relation in \eqref{energyrelationkj}.  In fact, the structure  of the energies provides a good 
topological measure for the potential solution; however the energy balance is {\it lost} in (\ref{energyrelationkj}) when making use of the state variable $\psi$ and, in addition, the boundary term involves the traces of $ L_2$ solutions of the flow, which are possibly {\it not defined}  at all.
In view of these complications,  our approach is be based on  (i) developing a suitable   theory for the traces of the flow solutions (as in \cite{supersonic}); (ii)  counteracting  the loss of energy balance relation. 

In this section we attempt to recover the dynamics after the change of state variable by viewing them as the sum of a generating component and a perturbation, where there perturbation corresponds to the ``energy polluting" term in the energy relation. This method was successfully used in \cite{supersonic}.

As explained in Section \ref{config}, the K-J boundary conditions arise in the modeling of panels that may be partially free. The existence of a pressure jump off of the wing occurs in many typical configurations \cite{bal0}. From the mathematical point of view
the difficulty lies, again, at the level of the linear theory. In order to deal with  the effects of the {\it unbounded}
traces $\gamma[\psi] $ in the energy relation \eqref{energyrelationkj} microlocal calculus is necessary. This
has been successfully accomplished in \cite{supersonic} where clamped boundary conditions in the supersonic case
were considered. However, in the case of K-J  boundary conditions there  is an additional difficulty that involves
``invertibility" of  {\it finite} Hilbert (resp. Riesz) transforms.  This latter property is known to fail within the $L_2$ framework, thus it is necessary to build the  $L_p$ theory, $p \ne 2 $.  This was for the first time observed in \cite{bal4}
and successfully resolved in the  one dimensional case. However, any progress to higher dimensions depends
on the validity of the corresponding harmonic analysis result developed for finite Hilbert/Riesz type transforms in two dimensions. The outline of a well-posedness proof are now outlined below.
\subsection{Abstract Setup}
We
 introduce the standard linear plate operator with {\em clamped boundary conditions}:
 $\mathscr{A}=\Delta^2$ with the domain
 \[
 \cD(\mathscr{A})=\{u\in H^4(\Omega): u=\Dn u = 0 \text{ on } \partial \Omega\}=(H^4\cap H_0^2)(\Omega).
   \]
  Additionally, $\mathscr{D}(\mathscr{A}^{1/2})= H_0^2(\Omega)
$. Take our state variable to be $$y\equiv(\phi, \psi; u, v) \in \big( H^1 (\realsthree_+) \times L_2(\realsthree_+)\big)\times\big( \cD(\cA^{1/2})\times L_2(\Omega)\big)\equiv Y.$$ We work with $\psi$ as an independent state variable, i.e., we are not explicitly taking $\psi = \phi_t+U\phi_x$ here. 

To build our abstract model, we must first define a preliminary operator $\bA: \cD(\bA) \subset Y \to Y$ by
\begin{equation}\label{op-A}
\bA \begin{pmatrix}
\phi\\
\psi\\ u\\ v
\end{pmatrix}=\begin{pmatrix}-U\partial_x \phi+\psi\\- U\partial_x\psi + \Delta \phi\\ v \\ -\cA u+ \gamma[\psi] \end{pmatrix}
\end{equation}
The domain of $\cD(\bA) $ is given by

\begin{equation}\label{dom-bA}
\cD(\bA) \equiv \left\{ y =  \begin{pmatrix}
\phi\\
\psi\\ u\\ v
\end{pmatrix}\in Y\; \left| \begin{array}{l}
-U \Dx \phi + \psi \in H^1(\R^3_+),~\\ -U \Dx \psi  + \Delta \phi  \in L_2(\R^3_+), \\
\psi =0 \text{ on } \realstwo\backslash\Omega, ~~\Dn \phi =- v , \text{ in} ~ \Omega \\
v \in \cD(\cA^{1/2} )= H_0^2(\Omega),~\\
-\cA u + \gamma  [\psi] \in L_2(\Omega) \end{array} \right. \right\}
\end{equation}
The operator $\bA$ will be the foundation of our abstract setup, and ultimately the dynamics of the evolution in \eqref{flowplate2} can be represented through $\bA$. 

As a first step, we show that  $\bA$ is m-dissipative.  
 \subsection{Semigroup Generation of $\bA$}
\begin{theorem}
The operator $\bA$ is  $m$-dissipative on $Y$.  Hence, via the Lumer-Philips theorem, it generates a C$_0$ semigroup of contractions.
\end{theorem}
\begin{proof} 
\noindent{\it Maximal Dissipativity}.
We employ the following inner product on our state space: for $y,~\hat y \in Y$

$$ (( y, \hat{y}))_Y \equiv  ( \nabla \phi , \nabla \hat{\phi} )_{  \realsthree_+}  + (\psi, \hat{\psi})_{  \realsthree_+}   + < \cA^{1/2} u, \cA^{1/2} \hat{v} >_{\Om}  + < v, \hat{v} >_{\Om} $$

With this inner product, we dissipativity is checked in a straightforward manner utilizing the K-J boundary conditions and the coupling, and for any $y \in Y$ $$((\bA y,y))_Y=0.$$

We must verify the range condition for $\bA$.  
We will take the system corresponding to $\bA+\lambda I$, with $\lambda = 0$ first.
Given $(f_1,f_2;g_1,g_2) \in Y$, we consider:
\begin{eqnarray}
-U \phi_x + \psi = f_1 \in H^1(  \realsthree_+) \\
 -U \psi_x + \Delta \phi = f_2\in L_2(  \realsthree_+) \\
v = g_1 \in \cD(\cA^{1/2} ) \\
- \cA u + \gamma[\psi] = g_2 \in L_2(\Om),
\end{eqnarray}
with boundary conditions:
\begin{equation}\begin{cases}\Dn \phi = -v = g_1 \in H^2(\Omega)\\ 
\psi =0 \text{ in }~ \realstwo\backslash\Omega.\end{cases}\end{equation}

Denoting $\kappa \equiv \phi_x $ and $\Delta_U \equiv \Delta - U^2 D_x^2 $  we easily deduce that $\kappa$ satisfies the following Zaremba (mixed) problem: 

$$\Delta_U \kappa = f_{2x} + U f_{1xx}\in H^{-1}(R^3_{+}) , ~~
\Dn \kappa =-v_x = g_{1x} \in H^1(\Omega) , ~~
 \kappa =  -\frac{1}{U}  f_1 \in H^{1/2} (\realstwo\setminus \Om),$$

As we are in the subsonic case, $ U < 1$, the above problem defined by $ \Delta_U $  is strongly elliptic in $ \kappa=\phi_x$. 

This is a mixed (Zaremba) problem and we recover $\phi_x \in H^1(  \realsthree_+ ). $ Returning to the equation, we recover $\psi$: $ \psi = U \phi_x + f_1 \in H^1(  \realsthree_+) $ and hence,
 $\psi \in H^1(\realsthree_+) $. With $\psi\in H^1(\realsthree_+)$ in hand (and hence $\gamma[\psi] \in H^{1/2}(\Omega)$), solving for $(u,v)$ is standard. In addition, having solved for $\phi_x$, we may then specify that $\Delta \phi = f_2+\psi_x \in L_2(\realsthree_+)$, with appropriate boundary conditions. We must verify that this is valid, by recovering $\phi \in H^1(\realsthree_+)$.
 
\begin{remark} Note that (from the regularity of the flow equations and mixed boundary conditions) we {\em do not obtain} $\phi \in H^2(\realsthree_+)$, demonstrating that the resolvent operator is not compact. \end{remark}
 
 To see that $\phi \in H^1 (\realsthree_+) $, we proceed as follows: let $\lambda > 0 $ and consider the equation for $(\bA +\lambda I)y=(f_1,f_2;g_1,g_2) \in Y$, where $y=(\phi,\psi;u,v)$ is a solution (as obtained above for the case $\lambda=0$):
 \begin{align*}
-U \phi_x + \psi  +\lambda \phi =& f_1 \in H^1(  \realsthree_+) \\
 -U \psi_x + \Delta \phi  + \lambda \psi =& f_2\in L_2(  \realsthree_+) \\
v  + \lambda u =&g_1 \in \cD(\cA^{1/2} ) \\
- \cA u + \gamma[\psi] + \lambda v  =& g_2 \in L_2(\Om)\\
\Dn \phi = -v , ~\text{in}~ \Omega , ~~\psi =&0 \text{ in }~ \realstwo\backslash\Omega,
\end{align*}
where $\cA$ stands for plate generator.  
Now, using the same idea as in \cite{supersonic}, we may obtain the a priori estimate on the solution:
$$\lambda ||(\phi,\psi;u,v)||_Y^2 \le C||(f_1,f_2;g_1,g_2)||_Y^2.$$
\end{proof}
 In addition, we note from the proof of the $m$-dissipativity of $\bA$ above, that $-\bA$ is also m-dissipative; indeed, $((-\bA y,y))_Y=0$, and the proof of maximality (the corresponding estimates) does not depend on the sign of $\bA$, owing to inherent cancellations in the structure of the static flow problem. Thus, with both $\pm \bA$ m-dissipative, we have: \begin{corollary}
The operator $\bA$ is skew-adjoint on $Y$ and generates a C$_0$ group of isometries.
 \end{corollary}

 \subsection{Encoding the Flow Boundary Conditions}
  In order to incoporate the flow boundary conditions abstractly into our operator representation of the evolution, we introduce the {\em flow-Neumann map}
 defined for the flow operator $$\bA_0 \begin{pmatrix}\phi\\\psi\end{pmatrix} \equiv   \begin{pmatrix}- U \phi_x + \psi\\ - U \psi_x + \Delta \phi \end{pmatrix} , $$ 
 with 
 $$\cD(\bA_0)= \{ (\phi, \psi) \in Y_f \equiv H^1(  \realsthree_+) \times L_2(  \realsthree_+)\big|~ 
 - U \phi_x + \psi \in H^1(  \realsthree_+),  ~U \psi_x + \Delta \phi\in L_2 (  \realsthree_+), $$ $$\Dn \phi =0 \text{ on }~ \Omega, ~\psi =0  ~\text{ in }~\realstwo\backslash\Omega \}.$$
By the  utilizing the arguments above in the proof of maximality, we have that  the operators $\pm\bA_0$ are $m$-dissipative on $Y_f$. This indicates that $\bA_0$ is skew-adjoint and one may check that 
 $$ (( \bA_0 y, \hat{y}))_{Y_f} =- ((y , \bA_0 \hat{y} ))_{Y_f}  ~\text{ for }~ y, ~\hat{y} \in \cD(\bA_0).$$

With a help of the flow map $\bA_0$ discussed above,  we define Neumann-flow map as follows:
$$ N : L_2(\realstwo) \rightarrow Y_f $$ given by 
$(\phi, \psi ) = Ng , ~~\text{ iff }$
$-U \phi_x + \psi + \phi  =0 ~\text{and} ~- U \psi_x + \Delta \phi + \psi  =0~ \text{ in }~   \realsthree_+, $ \noindent with $
\psi =0 ~\text{ in }~\realstwo\backslash\Omega ~\text{and}~\Dn \phi =-g  ~\text{ in }~ \Omega.$

We then consider the associated regularity of the map $N$. Note that the Neumann map  is associated with the matrix operator $\bA_0 $ rather than the usual harmonic extensions associated with a scalar elliptic operator. This difference is due to the fact that K-J conditions affect both the flow and the  aeroelastic potential. 
In order to describe the regularity of $N$ map we shall use the following anisotropic function spaces: 
$$H^{r,s}(D) = \{ f \in H^s(D),  \frac{\partial  f^r}{\partial x^r } \in H^s(D) \} $$ 
These spaces are subspaces of $H^s(D) $ with the additional information on regularity in $x$-direction. Using the regularity associated to the Zaremba elliptic problem (as discussed above) we have:
\begin{lemma}
$$N\in \mathscr{L} \big(H^{1,-1/2}(\Omega) \rightarrow H^{1,1}(  \realsthree_+) \times H^{0,1}(  \realsthree_+)\big) $$ 
where $g\in H^{1,-1/2}(\Omega) $ denotes the anisotropic space $\partial_xg \in H^{-1/2}(\Omega)$ , $g \in H^1(\Omega) $.
\end{lemma}
 Our next result identifies $N^*[ \bA_0^*+I] $ with a trace operator.  Here, the adjoint taken is with respect to the $L_2(\Omega) \rightarrow Y_f$ topology. This is reminiscent of a classical Neumann map: 
 \begin{lemma}\label{duality}
 Let $(\phi, \psi) \in \cD(\bA_0^*) $. Then 
 $N^*[ \bA_0^* + I ]  (\phi, \psi) = \gamma [\psi] $.  
 \end{lemma}

 \begin{proof}
 This follows from the identification of $\cD(\bA_0) = \cD(\bA_0^*),$ (from the skew-adjointness property) and Green's theorem.

\end{proof}
With the introduced notation we can express the final flow-structure operator as 
 \begin{equation}\label{op-AS}
\bA \begin{pmatrix}
\phi\\
\psi\\ u\\ v
\end{pmatrix}=\begin{pmatrix} \cdot\\\bA_0  \Big[\begin{pmatrix}\phi\\ \psi\end{pmatrix}  - N v \Big]  - Nv 
 \\ v \\ -\cA u +N^* (\bA_0^*  + I )  \begin{pmatrix}\phi\\\psi \end{pmatrix} \end{pmatrix}.
\end{equation}
This new representation of $\bA$ encodes the boundary conditions, and further reveals the antisymmetric structure of the problem.

\subsection{Cauchy Problem}
Having established that $\bA$ is $m$-dissipative, the Cauchy problem
\begin{equation}\label{eq-trankt-0}
y_t=\bA y, ~y(0)=y_0 \in Y
\end{equation}
 is well-posed on $Y$. However, this semigroup statement is not equivalent to the problem in \eqref{flowplate2}. The dynamics of the {\it original } fluid-structure interaction in \eqref{flowplate2} can be re-written (taking into account the action and domain of $\bA$) as
\begin{equation}\label{op-ASf}
\bB_L   \begin{pmatrix}
\phi\\
\psi\\ u\\ v
\end{pmatrix}  = \bA \begin{pmatrix}
\phi\\
\psi\\ u\\ v
\end{pmatrix}  + \bP \begin{pmatrix}  \phi \\ \psi\\ u\\ v \end{pmatrix} ,
\end{equation}
where $\bP$ corresponds to what remains of the dynamics in \eqref{flowplate2} which are not captured by $\bA$.
This allows us to treat the problem of well-posedness within the framework of ``unbounded trace perturbations",
where the perturbation in question becomes 
\begin{equation}\label{op-ASP}
\bP \begin{pmatrix}
\phi\\
\psi\\ u\\ v
\end{pmatrix}=\begin{pmatrix} \cdot\\-U (\bA_0+I)N u_x 
 \\ 0 \\ 0  \end{pmatrix}
\end{equation}
Here $(\bA_0+I)N $ is defined via duality (using its adjoint expression) via Lemma \ref{duality}.

We now verify that $\bA+\bP$ (computed formally) fully encodes the dynamics of \eqref{flowplate2}. 
\begin{align*}
(\bA+\bP)y=&\begin{pmatrix}\cdot\\ \bA_0\Big[\begin{pmatrix} \phi \\ \psi\end{pmatrix}-N(v+Uu_x)\Big]-N(v+Uu_x) \\ v \\ -\cA+N^*(\bA^*_0+I)\begin{pmatrix}\phi\\\psi \end{pmatrix} \end{pmatrix}
\end{align*}
That the plate components are correct is standard. We focus on the flow component: 
let $\begin{pmatrix} \hat \phi \\ \hat \psi \end{pmatrix} = -N(v+Uu_x)$. This implies that 
\begin{equation}\label{piece}-U\hat \phi_x+\hat \psi = - \hat\phi,~~-U\hat\psi_x+\Delta \hat\phi = -\hat\psi, ~~\Dn \hat\phi = -(v+Uu_x)~\text{ on } \Omega,~~\hat\psi = 0 \text{ on } \realstwo\backslash \Omega.\end{equation}
Hence
\begin{align*}
\bA_0\Big[\begin{pmatrix} \phi \\ \psi\end{pmatrix}-N(v+Uu_x)\Big]-N(v+Uu_x) =
& ~\bA_0\Big[ \begin{pmatrix} \phi \\ \psi \end{pmatrix}+\begin{pmatrix} \hat \phi \\ \hat \psi\end{pmatrix}\Big]
+\begin{pmatrix} \hat \phi \\ \hat \psi\end{pmatrix}
\\
=&~\begin{pmatrix}-U\phi_x+\psi \\ -U\psi_x+\Delta \phi \end{pmatrix} 
\end{align*} where we have used \eqref{piece} in the last line to make the cancellation.

\subsection{Generation for the Full Dynamics: $\bA+\bP$}
We would like to recast the full dynamics of the problem in \eqref{flowplate2} as a Cauchy problem in terms of the operator $\bA$. To do this, we define an operator $\bP:Y \to \mathscr{R}(\bP)$ as follows:
\begin{equation}\label{op-P}
\bP\begin{pmatrix}\phi\\\psi\\u\\v \end{pmatrix}= \bP_{\#}[u]\equiv\begin{pmatrix}0\\-U\bA_0N\partial_x u\\0\\0\end{pmatrix}
\end{equation}
Specifically, the problem in (\ref{flowplate2}) has the abstract Cauchy formulation:
\begin{equation*}
y_t = (\bA +\bP)y, ~y(0)=y_0,
\end{equation*}
where $y_0 \in Y$ will produce semigroup (mild) solutions to the corresponding integral equation, and $y_0 \in \cD(\bA)$ will produce classical solutions. To find solutions to this problem, we will consider a fixed point argument, which necessitates interpreting and solving the following inhomogeneous problem,  and then producing the corresponding estimate on the solution:
\begin{equation} \label{inhomcauchy}
y_t = \bA y +\bP_\# \overline{u}, ~t>0, ~y(0)=y_0,
\end{equation}
for a given $\overline{u}$. To do so, we must understand how $\bP$ acts on $Y$
(and thus $\bP_\#$ on $H^2(\Om)$).
\par
To motivate the following discussion, consider for $y \in Y$ and $z= (\overline\phi, \overline\psi;  \overline u, \overline v)$  and apply Lemma \ref{duality}  (with $Y$ as the pivot space) to obtain:
\begin{align}\label{p-y-z}
(\bP y, z)_Y =& (\bP_{\#}[u],z)_Y
= -U(\bA_0N\partial_x u,\overline \psi)
=-U<\partial_x u, \gamma[\overline \psi]>.
\end{align}
Hence, interpreting the operator $\bP$ (via duality) is contingent upon the ability to make sense of $\gamma[\overline \psi]$, which \textit{can} be done if $\gamma[\overline \psi] \in H^{-1+\epsilon}(\Omega)$ (since $u_x \in H^1(\Omega)$). In what follows, we show a trace estimate on $\psi$ (for semigroup solutions) of (\ref{inhomcauchy}) allows us to justify the program outlined above. We note that similar arguments were used in \cite{supersonic}. 
\par
We now state the trace regularity which is required for us to continue the abstract analysis of the dynamics. We state the following as a theorem which depends upon an assumption about the integral transform on $\Omega$. We will then show that abstract assumption is satisfied when the problem is reduced to a one dimensional structure.
  \begin{theorem}[{\bf Flow Trace Regularity}]\label{le:FTR0}
Let the auxiliary Assumption \ref{hilb2} stated below be in force. If $\phi(\xb,t)$  satisfies \eqref{flow} taken with flow boundary conditions \eqref{KJ}, then for every $\epsilon > 0 $ 
with $\psi = \phi_t + U \phi_x $ we have
\begin{equation}\label{trace-reg-est-M0}
\int_0^T\|\gamma[\psi](t)\|^2_{H^{-1/2 -\epsilon} (\realstwo)}dt\le C_T\left(
E_{fl}(0)+
 \int_0^T\| \Dn \phi(t) \|_{\Omega}^2dt\right)
\end{equation}
\end{theorem}
In Section \ref{tracereg} we state and give an explanation of Assumption \ref{hilb2}. 
\begin{remark}
We note that similar trace result has been established for the model considered in   \cite{supersonic} where, however, 
in this treatement the authors were able to take $\epsilon =0$. Since the above trace result is used in order to justify the  duality pairing  in (\ref{p-y-z}), 
it suffices to take any $\epsilon < 1$. 
\end{remark}

To conclude the proof of generation, we utilize the approach taken in \cite{supersonic} by interpreting the variation of parameters formula for $\overline{u} \in C(\R_+; H^2_0(\Om))$
\begin{equation}\label{solution1}
y(t)=e^{\bA t}y_0 + \int_0^t e^{\bA(t-s)}\bP_{\#} [\overline{u}(s)] ds.
\end{equation}
by writing (with some $\la\in\R$, $\la\neq 0$):
\begin{equation}\label{solution2}
y(t)=e^{\bA t}y_0 + (\lambda  -\bA)\int_0^t e^{\bA(t-s)}(\lambda  - \bA)^{-1}\bP_{\#} [\overline{u}(s)]  ds.
\end{equation}
 We initially take this solution in $[\cD(\bA^*)]'=[\cD(\bA)]'$ (via the skew-adjointness of $\bA$), i.e., by considering the solution $y(t)$ in (\ref{solution2}) above acting on an element of $\cD(\bA^*)$.

The approach is based on a key theorem in the theory of abstract boundary control \cite[pp.645-653]{redbook}. This theorem allows us to view the operator $\bP$ (mapping $Y$ outside of itself) as an unbounded perturbation. To make use of it, we critically implement the trace regularity theorem above in \eqref{le:FTR0}.  
We may now consider mild solutions to the problem given in (\ref{inhomcauchy}).
Applying general results on $C_0$-semigroups we arrive at the following assertion.
\begin{proposition}\label{pr:mild}
Let $\overline{u} \in C^1([0,T];H^{2}_0(\Omega))$ and $y_0\in Y$. Then $y(t)$
given by \eqref{solution1} belongs to  $C([0,T];Y)$ and
is a strong solution to \eqref{inhomcauchy} in $[\cD(\bA)]'$, i.e.
in addition we have that
\[
y\in  C^1((0,T);[\cD(\bA)]')
\]
and  \eqref{inhomcauchy} holds in $[\cD(\bA)]'$ for each $t\in (0,T)$.
\end{proposition}

From here, we utilize the full strength of the abstract semigroup theory with unbounded trace perturbation \cite[pp.645-653]{redbook} (and restated and used in \cite{supersonic}). 
\begin{theorem}[{\bf $L$ Regularity}]
Let $T>0$ be fixed, $y_0 \in Y$ and $\overline u \in C([0,T]; H_0^{2}(\Omega))$.
Then  the mild solution
 \begin{equation*}
  \ds y(t)=e^{\bA t}y_0+L[\overline u](t)\equiv
 e^{\bA t}y_0+\int_0^te^{\bA (t-s)}\bP_{\#}[\overline u(s)] ds
 \end{equation*}
 to problem \eqref{inhomcauchy} in $[\cD(\bA)]'$ belongs to the class $C([0,T]; Y)$
 and enjoys the estimate
\begin{align}\label{followest}
\max_{\tau \in [0,t]} ||y(\tau)||_Y \le&  ||y_0||_Y + k_T ||\overline u||_{L_2(0,t;H^{2}_0(\Omega))},~~~\forall\, t\in [0,T].
\end{align}
\end{theorem}
\begin{remark}
We emphasize that the perturbation $\bP$ acting outside of $Y$ is regularized when incorporated into the operator $L$ defined above; namely, the variation of parameters operator $L$ is a priori only continuous from $L_2(0,T;\cU)$ to $C(0,T;[\cD(\bA ^*)]')$. However, we have shown that the additional ``hidden'' regularity of the trace of $\psi$ for solutions to (\ref{flow}) with the boundary conditions in \eqref{KJ} allows us to bootstrap $L$ to be continuous from $L_2(0,T;\cU)$ to $C(0,T;Y)$ (with corresponding estimate) via the abstract theorem \cite{redbook}. This result  justifies \textit{formal} energy methods on the equation (\ref{inhomcauchy}) in order to produce a fixed point argument---see \cite{supersonic}.\end{remark}

\subsection{Construction of a Generator}

Let  $\bX_t = C\big((0,t];Y\big).$
Now,  with $y_0 \in Y $  take  $\overline{y}=(\overline{\phi},\overline{\psi};\overline{u},\overline{v}) \in \bX_t$, and introduce the map $\cF: \overline{y} \to y$ given by
\begin{equation*}
y(t) = e^{\bA t}y_0+L[\overline u](t),
\end{equation*}
i.e. $y$ solves 
$y_t=\bA y+\bP_\# \overline u,
~~ y(0)=y_0,$
in the generalized sense, where $\bP_\#$ is defined in \eqref{op-P}.
Using (\ref{followest}) we then implement a standard fixed point argument on $\bX_t$.

This implies that there exists a $t_*$such that on the interval $[0,t_*]$ the problem
$
y_t=\bA y+\bP y, ~~t>0,~
~~ y(0)=y_0,
$
has a local in time unique (mild) solution defined now in $Y$. This above local solution
can be extended to a global solution in finitely many steps by linearity.
Thus there exists a unique function
$y=(\phi,\psi;u,v)\in C\big(\R_+;Y\big)$ which satisfies the variation of parameters (Duhamel) formula.
It also follows from the analysis above that
\[
\|y(t)\|_Y\le C_T \|y_0\|_Y,~~~ t\in [0,T],~~\forall\, T>0.
\]
Thus we have a strongly continuous semigroup corresponding to solutions
$\widehat{T}(t)$ in $Y$.
Additionally, since mild solutions satisfy the corresponding variational equality we have
\begin{equation*}
    (y(t),h)_Y=(y_0,h)_Y+ \int_0^t\left[
    -(y(\tau), \bA h)_Y+ (\bP[y(\tau)],h)_Y
    \right]d\tau
    ,~~~\forall\, h\in \cD(\bA),
    ~~ t>0.
\end{equation*}
Using the same idea as presented in \cite{jadea12, supersonic}, we conclude that
the generator $\widehat{\bA}$ of $\widehat{T}(t)$  has the form
\[
\widehat{\bA}z=\bA z+\bP z,~~z\in\cD(\widehat{\bA})=\left\{z\in Y\,:\; \bA z+\bP z\in Y\right\}
\]
(we note that  the sum $\bA z+\bP z$ is well-defined as an element
in $[\cD(\bA)]'$ for every $z\in Y$). Hence, the semigroup $e^{\widehat \bA t}y_0$ is a generalized solution for $y_0 \in Y$ (resp. a classical solution for $y_0 \in \cD(\widehat \bA)$) to (\ref{flowplate2}) on $[0,T]$ for all $T>0$.
\begin{equation}\label{dom-bA-n}
\cD(\bA + \bP) \equiv \left\{ y \in Y\; \left| \begin{array}{l}
-U \Dx \phi + \psi \in H^1(\R^3_+),~\\ -U \Dx \psi  -\bA_0 (\phi - N (v + U \Dx u )) \in L_2(\R^3_+) \\
v \in \cD(\cA^{1/2} )= H_0^2(\Omega),~
-\cA u +N^* \bA_0^* \psi \in L_2(\Omega) \end{array} \right. \right\}
\end{equation}

Indeed, the
function $y(t)$ is a generalized solution corresponding to the generator $\bA + \bP$ with the  domain defined in (\ref{dom-bA-n}).

\subsection{Trace Regularity}\label{tracereg}
As we noted above, the ``hidden" trace regularity of the term $\psi$ coming from the flow equation \eqref{flow} is {\em critical} to the arguments above. In this section we analyze this problem in the dual (Fourier-Laplace) domain and relate it to a certain class of integral transforms reminiscent of the finite Hilbert transform. In the case of two dimensions, we reduce the trace regularity to an hypothesis about the invertibility of Hilbert-like transforms on bounded domains. Additionally, we assert the necessary trace regularity on a pseduodifferential operator corresponding to the flow problem in one dimension. 

We are interested in the trace regularity of the following flow problem in $\realsthree_+$: \begin{equation}\label{flow1}\begin{cases}
(\partial_t+U\partial_x)^2\phi=\Delta \phi & \text { in } \realsthree_+ \times (0,T),\\
\phi(0)=\phi_0;~~\phi_t(0)=\phi_1,\\
\Dn \phi = d(\xb,t)& \text{ on } \Omega \times (0,T)\\
\gamma[\phi_t+U\phi_x]=\gamma[\psi]=0, ~&\xb \in \realstwo \backslash \Omega,
\end{cases}
\end{equation}
with $0 \le U <1$,  $d(\xb,t) $ is the ``downwash" generated on the structure, and $\gamma [\psi] $ is the aeroelastic potential. Aeroelastic potential satisfies the trace estimate given in Theorem \ref{le:FTR0}.
\begin{remark}
A related analysis of  trace regularity for this equation was carried out in the case 
of $L_2^{loc}(\realstwo)$ purely Neumann data in \cite{supersonic}. Here, the mixed boundary conditions present a  challenge in the microlocal analysis.
\end{remark}

The analysis of the trace regularity is done in the dual (Fourier-Laplace) domain. It is here where we encounter the analysis (microlocally) involving the an integral transform which is analogous to the finite Hilbert transform on $\Omega$. We now introduce the necessary assumption mentioned above. Define and operator $\mathscr{H}$ on $L_2(\Omega)$ whose symbol is given by $symb(\mathscr H)= \dfrac{-i|\eta_x|}{\eta_x}\equiv j(\eta_x)$, where $\eta_x \leftrightarrow \frac{1}{i}\partial_x$ in the Fourier correspondence. Additionally, let $P_{\Omega}:\reals^2 \to \Omega$ be the associated projection into $L_p(\Omega)$ and $E_{\Omega}$ be the extension operator (by zero) from $L_p(\Omega)$ into  $L_p(\realstwo) $.   

\begin{assumption}\label{hilb2}
Assume that the operator $\mathscr{H}_f\equiv P_{\Omega}\mathscr{H}E_{\Omega}:L_p(\Omega) \to L_p(\Omega)$ is continuously invertible for $p \in (1,2)$. 
\end{assumption}
When $\Omega\equiv (-1,1)$ and the problem is reduced to a two dimensional flow interacting with the one dimensional structure, we have that Assumption \ref{hilb2} reduces to the invertibility of the finite Hilbert transform \cite{FH1,FH,FH2}. That the finite Hilbert transform is invertible on $L_{2^-}(-1,1)$ is discussed in \cite{FH1,FH}. We note that it is also Fredholm on $L_{p}(-1,1)$ for $p>2$, and in the case $p=2$ the finite Hilbert transform has a range which is proper and dense. Hence, when $\Omega = (-1,1)$, Assumption \ref{hilb2} is satisfied. 

 \begin{remark}
The connection between integral equations appearing in the study of aeroelasticity and invertibility of finite Hilbert fransforms has been known for many years and dates back to Tricomi and his airfoil equation \cite{tricomi}. This approach has been critically used in \cite{bal0,bal4} where the analysis is centered on solvability of integral equations connecting the downwash with the aeroelastic potential. In these works the author performs an analysis on a one dimensional flow-beam system and utilizes a similar Fourier-Laplace approach and corresponding $L_p$ theory. However, the solution to the system (the Possio integral equation) is {\em constructed} via a Fourier-Laplace approach. We follow the same conceptual route with different technical tools. Our approach is based on microlocal analysis, rather than explicit solvers of integral equations arising in a purely one dimensional setting. Though our final estimate depends on an assumption (which we demonstrate is satisfied only for the one dimensional 
case) we 
believe that the microlocal approach provides new ground for extending the flow-structure analysis to the multidimensional settings.
 \end{remark}
 We note at this stage that the integral transform which arises corresponding to the so called  abstract version of Possio equation \cite{bal0,bal4}  $$d= P_{\Omega} \mathscr{H} \mathscr{S}
  E_{\Omega} f; x\in \Omega , t > 0 $$ where $f$ is supported only on $\Omega$. In our case, via the K-J condition, $\psi$ is a function on the whole domain whose support lies inside $\Omega$; hence we have:
$$d=P_{\Omega}\mathscr H \mathscr S E_{\Omega}  \psi,~\text{on}~ \Omega , t > 0 ,$$
We recall that    $\psi =0$ off $\Omega$, so $E_{\Omega}\psi $ coincides in this case with $\psi$.
Given $d \in C(0, T; L_2(\Omega) ) $ our task is to infer the regularity of $E_{\Omega} \psi $. 
This precise formulation is the one given in \cite{bal0}, where the connection with finite Hilbert transform has been used. 
 This is {\em not a standard  problem in the two dimensional scenario}. However, when $\Omega$ is reduced to an interval, the integral equation above contains the classical finite Hilbert transform. 
 The pseudodifferential operator $\mathscr{S}$ (corresponding to an auxiliary symbol $S(\beta,\eta)$---$\beta$ being dual time variable and $\eta$ being the spatial dual variable) is an operator with the property that 
 $$ \mathscr{S}^{-1}: L_2(0,T; H^{-\epsilon}  (\reals) )   \rightarrow L_2(0,T; H^{-1/2 - \epsilon}(\reals))  , \epsilon > 0  $$  is bounded (a familiar loss of $1/2 $ derivative \cite{supersonic}).
Relating the Hilbert transform $\mathscr H$ to the finite Hilbert transform described above one obtains 
\begin{align}\label{possio}
d =  P_{\Omega} \mathscr{H} \mathscr{S} E_{\Omega} \psi = P_{\Omega} \mathscr{H} E_{\Omega} P_{\Omega } \mathscr{S} E_{\Omega} 
   \psi \nonumber \\
 +P_{\Omega} \mathscr{H} [ I - E_{\Omega} P_{\Omega } ]  \mathscr{S} E_{\Omega} \psi \nonumber \\
 = \mathscr{H}_f P_{\Omega} \mathscr{S} E_{\Omega} \psi + V  \mathscr{S} E_{\Omega} \psi 
 \end{align}
 Since the singular support of $ [ I - E_{\Omega} P_{\Omega} ] \mathscr{S} E_{\Omega} \psi $ is empty,
 by the pseudolocal properties of pseudodifferential operators,  the operator $V$ is ``smooth " and compact (see  \cite{bal0,bal4} for detailed calculations with a similar decomposition). The operator $\mathscr{H}_f P_{\Omega}   + V $  is then invertible on  $L_2(0, T; L_p(\Omega) ) $ with $ p < 2$. 
  Thus  with $d \in L_2(0,T; L_2(\Omega))$, by the properties of finite Hilbert Transform,  we have $P_{\Omega}  \mathscr{S} \psi \in L_2(0, T; L_{p}(\reals))$ with $p < 2 $, which yields via Sobolev's embeddings $ \mathscr{S} \psi \in L_2(0, T;  
  H^{-\epsilon} (\Omega) )  $ for every $\epsilon> 0$ by taking a suitable $ p < 2 $. Thus $\psi \in L_2(0,T; H^{-1/2-\epsilon}(\Omega) ) $, as desired.

 A detailed discussion of the relationship between Theorem \ref{le:FTR0} and the integral transform $\mathscr{H}_f$ can be found in a forthcoming manuscript. This treatment will also contain the details of the microlocal proof outlined above, and a proof of the invertibility $\mathscr{S}$ when $\Omega=(-1,1)$.

\begin{remark}\label{Bal1}
We note that a similar result is obtained in the analysis in \cite{bal0,bal4}, where the author proves that 
aeroelastic potential $ \psi \in L_2(0, T;L_q(\Omega) ) $ for $ q < 4/3 $ . Since  for $p > 4 $ there exists $\epsilon > 0 $ 
such that $$H^{1/2 + \epsilon} (\Omega) 
\subset L_p(\Omega) , ~~p > 4 , ~~\text{dim}~~\Omega \leq 2,$$  and one then obtains  that $L_q(\Omega) \subset H^{-1/2 - \epsilon}(\Omega) $ with $q < 4/3 $. \end{remark}
\section{Long-time Behavior of Solutions}
\subsection{Reduction to a Delayed Model}
A key to obtaining the attracting property of the dynamics is the representation of the flow on the structure via a delay potential. We may rewrite the full flow-plate system as a von Karman system with delay term. Reducing the flow-plate problem to a delayed von Karman plate is the primary motivation for our main result and allows long-time behavior analysis of the flow-plate system.
The exact statement of this reduction is given in the following assertion:

\begin{theorem}\label{rewrite}
Let the hypotheses of Theorem~\ref{th:nonlin} be in force,
and $(u_0,u_1;\phi_0,\phi_1) \in \cH \times L_2(\Omega) \times H^1(\realsthree_+) \times L_2(\realsthree_+)$. 
Assume that there exists an $R$ such that $\phi_0(\xb) = \phi_1(\xb)=0$ for $|\xb|>R$.  
Then the there exists a time $t^{\#}(R,U,\Omega) > 0$ such that for all $t>t^{\#}$ the weak solution $u(t)$ to 
(\ref{flowplate}) in the clamped case  satisfies the following equation:
\begin{equation}\label{reducedplate}
u_{tt}+\Delta^2u-[u,v(u)+F_0]=p_0-(\partial_t+U\partial_x)u-q^u(t)
\end{equation}
with
\begin{equation}\label{potential}
q^u(t)=\dfrac{1}{2\pi}\int_0^{t^*}ds\int_0^{2\pi}d\theta [M^2_{\theta}\widehat u](x-(U+\sin \theta)s,y-s\cos \theta, t-s).
\end{equation}
Here, $\widehat u$ is the extension
 of $u$
   by 0 outside of $\Omega$; $M_{\theta} = \sin\theta\partial_x+\cos \theta \partial_y$ and \begin{equation*}
   t^*=\inf \{ t~:~\xb(U,\theta, s) \notin \Omega \text{ for all } \xb \in \Omega, ~\theta \in [0,2\pi], \text{ and } s>t\}
\end{equation*} with $\xb(U,\theta,s) = (x-(U+\sin \theta)s,y-s\cos\theta) \subset \realstwo$.
\end{theorem}
\begin{remark}
This extremely helpful theorem first appeared as a heuristic in \cite{kras} and was used in this way for many years; it was later made rigorous in \cite{Chu92b}.
\end{remark}

Thus, after some time, the behavior of the flow can be captured by the aerodynamical pressure term $p(t)$ in the form of a reduced delayed forcing.  Theorem~\ref{rewrite} allows us (assuming that the flow data is compactly supported) to suppress the dependence of the problem on the flow variable $\phi$.
Here we emphasize that the structure of aerodynamical pressure \eqref{aero-dyn-pr} posited in the hypotheses leads to the velocity term  $-u_t$ on the RHS of \eqref{reducedplate}.
We may utilize this as natural damping appearing in the structure of the reduced flow pressure by moving this term to the LHS.

\par
As we see below, the reduction method above allows us to study long-time behavior of the dynamical system corresponding to (\ref{flowplate}) (for sufficiently large times) by reducing the problem to a plate equation with delay. The flow state variables $(\phi,\phi_t)$  manifest themselves in our rewritten system via the delayed  character of the problem; they appear  in the initial data for the delayed  component of the plate, namely $u^t\big|_{(-t^*,0)}$.  Hence the behavior of both dynamical systems agree for all $t>t^{\#}(R,U,\Omega)$. By the dynamical systems property for solutions to the full system (semigroup well-posedness), we can propagate forward and simply study the long-time behavior of the plate with delay on the interval $(\sigma-t^*,\sigma+T]$ for $\sigma>t_{\#}$ and $T \le \infty$.
\begin{remark}
An immediate observation is that  both the nonlinear term $[u, v(u)+F_0]$  and the forcing term due to the flow,
$q^u$, are at the critical level with respect to the topology of phase space.
This immediately rules out the possibility of relying on {\em compactness} for the forcing terms---a critical property in studying long-time behavior.
In fact, coping with this issue presents the major challenge in solving the problem of long-time behavior of trajectories. 
\end{remark}
\smallskip\par
The following proposition concerns the delayed force term in the delayed von Karman plate  model (\ref{reducedplate}) and permits much of the analysis to follow.

\begin{proposition}\label{pr:q}
Let $q^u(t)$ be given by (\ref{potential}). Then \begin{equation}\label{qnegest}
||q^u(t)||^2_{-1} \le ct^*\int_{t-t^*}^t||u(\tau)||^2_1d\tau
\end{equation} for any $u \in L_2(t-t^*,t;H_0^1(\Omega))$.
If $u \in L_2^{loc}([-t^*,+\infty);H^2\cap H_0^1)(\Omega))$ we also have
\begin{equation}\label{qnegest2}
||q^u(t)||^2 \le ct^*\int_{t-t^*}^t||u(\tau)||^2_2d\tau,~~~\forall t\ge0,
\end{equation}
and
\begin{equation}\label{qnegest3}
\int_0^t ||q^u(\tau)||^2 d\tau \le c[t^*]^2\int_{-t^*}^t||u(\tau)||^2_2d\tau
,~~~\forall t\ge0.
\end{equation}
Moreover if $u \in C(-t^*,+\infty;H^2\cap H_0^1)(\Omega))$, we have that  $q^u(t) \in C^1( \R_+; H^{-1}(\Omega))$,
\begin{equation}\label{qnegest4}
\|q^u_t(t)\|_{-1} \le C\Big\{ ||u(t)||_1+||u(t-t^*)||_1+\int_{-t^*}^0||u(t+\tau)||_2d\tau\Big\},
~~~\forall t\ge0.
\end{equation}
\end{proposition}
\begin{proof}The proof of the bounds \eqref{qnegest}--\eqref{qnegest3} can be found in \cite{Chu92b} and \cite{springer}.
The estimate \eqref{qnegest4} follows by differentiating $q$ in time (distributionally). 
Details are presented in \cite{delay}. 
\end{proof}

\begin{remark}
{\rm
A priori, when $u_t$ is in $H^1_0(\Omega)$, it is clear from \eqref{qnegest} that there is a compactness margin  and
 we have the estimate
 $$
 \int_0^t<q^u(\tau),u_t(\tau)>d\tau\le
 \epsilon \int_0^t ||u_t(\tau)||^2_1+C(\epsilon,t)\sup_{\tau \in [-t^*,t]}||u(\tau)||^2_{1}.
 $$
 However, this is not immediately apparent when $u_t \in L_2(\Omega)$ as $||q^u(t)||^2_0$ has no such a priori bound from above, as in \eqref{qnegest}. This is precisely the property which disallows the previous ``abstract" analysis of second order equations with delay.
Hence, the critical component which will allow us to perform an analysis in the $\alpha =0$ case is the ``hidden compactness" of the term displayed in \eqref{qnegest4}; this allows estimation at the energy level.

  We note that inequality (\ref{qnegest4}) represents a loss of one derivative  (anisotropic---time derivatives are scaled by two spatial derivatives), versus the loss of two derivatives in  (\ref{qnegest}), (\ref{qnegest2}), and (\ref{qnegest3}).
}
\end{remark}

\subsection{Plate Model with Delay as a Dynamical System   }
Below we utilize a positive parameter $0<t^*<+\infty$ as the time of delay, and
accept  the commonly  used (see, e.g.,   \cite{Delay-book1995} or \cite{wu-1996}) notation $u^t(\cdot)$
 for function on $s\in [-t^*,0]$ of the form $s\mapsto u(t+s)$.
This is necessary due to the delayed character of the problem
which requires initial data on the prehistory interval  $[-t^*,0]$, i.e., we
need to impose an initial condition
 of the form $u|_{t \in (-t^*,0)} = \eta(\xb, t)$,
 where $\eta$ is a given function on $\Om\times [-t^*,0]$.
 We can choose this prehistory data $\eta$ in various
 classes. In our problem it is convenient to deal
 with Hilbert type structures, and therefore we assume in what follows that
 $ \eta \in L_2(-t^*,0;\cH)$. Since we do not assume
 the continuity of $\eta$ in $s\in [-t^*,0]$,  we also
 need to add the (standard) initial conditions of the form
 $u(t=0)= u_0(\xb)$ and $\partial_t u(t=0)=u_1(x)$.
 \par
Our delayed system is then given by:
\begin{equation}\label{plate}\begin{cases}
u_{tt}+\Delta^2u+u_t+f(u) +Lu= p_0+q(u^t,t) ~~ \text { in } ~\Omega\times (0,T), \\
u=\Dn u=0  ~~\text{ on } ~ \partial\Omega\times (0,T),  \\
u(0)=u_0,~~u_t(0)=u_1,~~\\ u|_{t \in (-t^*,0)} = \eta\in L_2(-t^*,0;H^2_0(\Omega)).
\end{cases}
\end{equation}
Here $f(u)$  is the von Karman nonlinearity as given above, and the forcing term $q(u^t,t)$ occurring on the RHS of the plate equation encompasses the gas flow (as in \eqref{rewrite}).
  The operator $L$ encompasses spatial lower order terms which do not have gradient structure
 (e.g., the term $-Uu_x$  in \eqref{reducedplate}). 
 
\subsubsection{Properties of the Delayed System}

 We take a \textit{weak solution} to \eqref{plate} on $[0,T]$ to be a function $$u \in L_{\infty}(0,T;H_0^2(\Omega))\cap W^1_{\infty}(0,T;L_2(\Omega)) \cap L_2(-t^*,0;H_0^2(\Omega))$$ such that the variational relation corresponding to \eqref{plate} holds (see, e.g., \cite[(4.1.39), p.211]{springer}).
We now assert:
\begin{proposition}\label{p:well}
The well-posedness of weak solutions to \eqref{flowplate}, as given in Theorem \ref{th:nonlin}, imply that the corresponding solutions are weak solutions to \eqref{plate}. Such solutions  to \eqref{plate} belong to
the class
$$C(0,T;\cH)\cap C^1(0,T;L_2(\Omega))
$$
and satisfies the energy identity
\begin{equation}\label{energyrelation2}
\cE(t)+\int_s^t ||u_t(\tau)||^2d\tau=\cE(s)+
\int_s^t<q(u^{\tau},\tau),u_t(\tau)>d\tau+\int_s^t<p_0-Lu(\tau),u_t(\tau)>d\tau,
\end{equation}
where the full (not necessarily positive) energy has the form
\begin{equation}\label{fulle}
\cE(u,u_t) \equiv \dfrac{1}{2}\big\{||u_t||^2-<[u,F_0],u>\big\}+\Pi_*(u)
\end{equation}
with \begin{equation}\label{potentiale}
\Pi_*(u) \equiv \dfrac{1}{2}\big\{ \|\Delta u\|^2+\dfrac{1}{2}||\Delta v(u)||^2\big\}.
\end{equation}
\end{proposition}
 \begin{lemma}\label{le:q} We denote $q^u(t)=q(u^t,t)$ and note the estimates in Proposition \ref{pr:q}. Then
\begin{align}\label{hidden1}
\Big|\int_0^t <q^u(\tau),u_t(\tau)> d\tau\Big| \le&~   C\e^{-1} t^*  \int_{-t^*}^t||u(\tau)||_2^2d\tau +\e \int_0^t ||u_t(\tau)||^2d\tau,~~~\forall \e>0,~\forall t\in[0,T],
\end{align}
for any  $u \in L_2(-t^*,T;H^2(\Omega))\cap W^1_2(0,T;L_2(\Omega))$.
\par
Additionally, there exists $\eta_*>0$ such that for every $\epsilon>0$  we have the estimate:
\begin{align}\label{hidden2} \Big|\int_0^t <q^u(\tau),u_t(\tau)> d\tau\Big| \le& ~ \epsilon\int_{-t^*}^t ||u(\tau)||_{2}^2d\tau +C(t^*,\e)\cdot(1+T)\sup_{[0,t]}||u(\tau)||^2_{2-\eta_*},~~\forall t\in[0,T],
\end{align}
 for  any $u \in L_2(-t^*,T;H^2(\Omega))\cap C(0,T;H^{2-\eta_*}(\Omega))\cap C^1(0,T;L_2(\Omega))$.
 \end{lemma}
 \begin{proof}
 The relation in \eqref{hidden1} easily follows from:
 \[
 \int_0^t d\tau \int_{\tau-t*}^\tau \phi(s) ds\le t^* \int_{-t*}^t \phi(s) ds,~~~\forall~\phi\in L_1(0,T).
 \]
 Inequality \eqref{hidden2} is verified by integrating by parts in $t$, applying \eqref{qnegest4} with $\psi = u(t)$,
and utilizing \eqref{qnegest2}.
For more details we refer to \cite{delay}.
\end{proof}

 In order to consider the delayed system as a dynamical system with the phase space $$\mathbf{H}\equiv \cH\times L_2(\Omega)\times L_2(-t^*,0;\cH),$$ we recall  the notation:
 $u^t(s) \equiv u(t+s) , s \in [-t^*, 0 ]$.
 With the above notation we introduce the  evolution operator
$ S_t\, : \bH\mapsto \bH$ by the formula
\begin{equation*}
  S_t(u_0, u_1, \eta) \equiv (u(t), u_t(t), u^t),
\end{equation*}
where $u(t)$ solves \eqref{plate}.
   Proposition \ref{p:well} implies the following conclusion
 \begin{corollary}\label{co:generation}
 $S_t : \mathbf{H} \rightarrow \mathbf{H}  $ is a strongly continuous semigroup on $\mathbf{H}$.
 \end{corollary}
\begin{proof}
Strong continuity is stated in Proposition~\ref{p:well}.
 The semigroup property  follows from uniqueness.
Continuity with respect to initial data {\em in this setting}  follows from:

\begin{lemma}[\cite{delay}]\label{l:lip}
Suppose $u^i(t)$ for $i=1,2$ are weak solutions to (\ref{plate}) with different initial data and $z=u^1-u^2$. Additionally assume that
\begin{equation}\label{bnd-R}
||u_t^i(t)||^2+||\Delta u^i(t)||^2 \le R^2, ~i=1,2
\end{equation}
 for some $R>0$ and all $t \in [0,T]$. Then there exists $C>0$ and $a_R\equiv a_R(t^*)>0$  such that
\begin{align*}
||z_t(t)||^2+||\Delta z(t)||^2 \le &~ Ce^{a_Rt}\Big\{||\Delta(u^1_0-u^2_0)||^2+||u^1_1-u^2_1||^2+\int_{-t^*}^0||\eta^1(\tau)-\eta^2(\tau)||_2^2d\tau\Big\}
\end{align*}
for all $t \in [0,T]$.
\end{lemma}
\begin{proof}
The argument is standard.  
We first note that $z$ solves the following problem
\begin{equation}\label{difference}
\begin{cases}z_{tt}+\Delta ^2 z + z_t+f(u^1)-f(u^2)=q(z^t,t)-Lz,
 \\ BC(z) \text{ on $\partial \Omega$ },
 \\ z(0)=z_0 \in \cH, ~z_t(0)=z_1 \in L_2(\Omega), ~z|_{(-t^*,0)}\in L_2(-t^*,0;\cH),
 \end{cases}
\end{equation}
and then apply the standard energy multiplier
and Gronwall's inequality.
\end{proof}
Using Lemma~\ref{l:lip} we obtain that
\begin{equation*}
||S_{t}y_{1}-S_{t}y_{2}||_{\bH}^{2}\leq C e^{a_Rt}||y_{1}-y_{2}||_{\bH}^{2}
\end{equation*}
for $S_ty_i=(u^i(t), u^i_t(t), [u^i]^t)$ with $y_i=(u^i_0, u^i_1, \eta)$,
such that \eqref{bnd-R} holds.
This allows us to conclude the proof of Corollary~\ref{co:generation}.
\end{proof}

 \subsection{Dissipative Dynamical System}
 In the case when the system is gradient system, existence of global attractors reduces to the study of asymptotic smoothness.  Since the model under consideration does not exhibit gradient structure (in the reduced, delayed setting), our first step in the  long-time analysis is showing an existence of an absorbing set; in other words, proving the dissipativity property.
Thus, our next task (in order to to make use of Theorem \ref{dissmooth}) is to show the dissipativity of the dynamical system $(S_t,\bH)$, namely that there exists a bounded, forward invariant, absorbing set. To show this, similar to the consideration in \cite[Theorem 9.3.4, p.480]{springer}, we consider the Lyapunov-type function (with $\cE(u,u_t)$ as in  (\ref{fulle}), and with $\Pi_*(u)$ given by \eqref{potentiale})
\begin{align*}
V(S_t y) \equiv &\cE(u(t),u_t(t))-<q^u(t),u(t)>+\nu\big(<u_t,u>+\frac{1+k}{2}||u||^2\big)\\\nonumber
&+\mu\Big(\int_{t-t^*}^t\Pi_*(u(s))ds+\int_0^{t^*}ds\int_{t-s}^t \Pi_*(u(\tau))d\tau\Big),
\end{align*}
where $S_t y \equiv y(t)= (u(t),u_t(t),u^t)$ for $t \ge 0$ and $\mu, \nu$ are some positive numbers to be specified below, and $k \ge 0$.
\begin{remark} The value $k$ above is the same as the (possible) damping coefficient added to 
the plate equation: $ku_t$. In what immediately follows we are not considering any 
imposed damping---the natural damping from the flow is being utilized via the reduction result---however, 
we list the above Lyapunov-type function in full generality. We restrict here to the case that $k=0$.
\end{remark}
In view of the results for the von Karman plate in \cite[Section 4.1.1]{springer}, we have that
\begin{equation}\label{energybounds}
c_0E(u,u_t) - c \le V(S_ty) \le c_1E(u,u_t)+\mu C t^*\int_{-t^*}^0 \Pi_*(u(t+\tau))d\tau +c
\end{equation}
for $\nu>0$ small enough,
 where $c_0,c_1,c,C >0$ are constant.
Here we make use  of the notation:
$
E(u,u_t) \equiv \dfrac{1}{2}||u_t||^2+\Pi_*(u).
$

 To obtain the above bound, we here need (and below) the critical lower bound on the potential energy cited in \eqref{potentialbound} (which can also be found on \cite[p. 49 and p. 132]{springer}). In what follows below, we will often make use of the above theorem to give $$<[u,F_0],u>+||u||^2 \le \delta \Pi_*(u)+C_{\delta,F_0}.$$
Next, the computation of $\ds \dfrac{d}{dt} V(S_ty)$ (presented in \cite{delay}) yields 

\begin{lemma}\label{le:48}
There exist $\mu, \nu >0$ and $c(\mu,\nu,t^*)>0$ and $C(\mu,\nu,p_0,F_0)>0$ such that
\begin{equation*}
\dfrac{d}{dt}V(S_ty)\le-c\big\{||u_t||^2+||\Delta u||^2+||\Delta v(u)||^2 +\Pi_*(u(t-t^*))+\int_{-t^*}^0\Pi_*(u(t+\tau))d\tau\big\}+C.
\end{equation*}
\end{lemma}

From this lemma and the upper bound in \eqref{energybounds}, we have for some $\beta>0$ sufficiently small (again, depending on $\mu$ and $\nu$): \begin{equation}\label{gronish}
\dfrac{d}{dt}V(S_ty) +\beta V(S_ty) \le C,~~t>0,
\end{equation}
The estimate above in (\ref{gronish}) implies (by a version of Gronwall's inequality) that
\begin{equation*}
V(S_ty) \le V(y)e^{-\beta t}+\dfrac{C}{\beta}(1-e^{-\beta t}).
\end{equation*}
Hence, the set
\begin{equation}\label{abs}
\mathscr{B} \equiv \left\{y \in \bH:~V(y) \le 1+\dfrac{C}{\beta} \right\},
\end{equation}  is a bounded forward invariant absorbing set.

Thus we have 
\begin{corollary}\label{c:a}
The system  $(\bH,S_t)$ is dissipative with absorbing set given by (\ref{abs})
\end{corollary}
We recall \cite{Babin-Vishik,ch-0,temam}
that
a closed set $B \subset \bH$ is said to be \textit{absorbing} for $(S_t,\bH)$ if for any bounded set $D \subset \bH$ 
there exists a $t_0(D)$ such that $S_tD \subset B$ for all $t > t_0$. 
If the dynamical system $(S_t,\bH)$ has a bounded absorbing set it is said to be \textit{dissipative}.
\par

\section{Proof of Theorem \ref{th:main2}}

In the context of this paper we will use a few keys theorems 
 to prove the existence of a finite dimensional global attractor.  First, we address
the existence of attractors and characterize the attracting set (for the proof, 
see  \cite{Babin-Vishik} or \cite{temam}):

\begin{theorem}
\label{dissmooth} Any dissipative and asymptotically smooth dynamical system $(S_t,\bH)$ in a Banach space $\bH$ possesses a unique compact global attractor $\textbf{A}$. This attractor is a connected set and can be described as a set of all bounded full trajectories.
\end{theorem}
We recall (see, e.g., \cite{Babin-Vishik,ch-0,lad,temam})   that
a \textit{global attractor} $\mathbf{A}$ is a closed, bounded set in
$\bH$ which is invariant (i.e., $S_t\mathbf{A}=\mathbf{A}$ for all
$t>0$) and uniformly attracts every bounded set $B$, i.e.
\begin{equation}  \label{dist-u}
\lim_{t\to+\infty}d_{\bH}\{S_t B|\mathbf A\}=0,~~~\mbox{where}~~
d_{\bH}\{S_t B|\mathbf A\}\equiv\sup_{y\in B}{\rm dist}_{\bH}(y,\mathbf A),
\end{equation}
for any bounded  $B\in\bH$.
One says that a dynamical system
$(S_t,\bH)$ is \textit{asymptotically smooth} if for any
bounded, forward invariant  set $D$ there exists a compact set $K
\subset \overline{D}$ such that
\begin{equation*} 
\lim_{t\to+\infty}d_{\bH}\{S_t D|K\}=0
\end{equation*}
 holds. An asymptotically smooth
dynamical system should be thought of as one which possesses \textit{local
attractors}, i.e. for a given forward invariant set $B_R$ of diameter $R$ in the space $
\bH$ there exists a compact attracting set in the closure of $B_R$,
however, this set need not be uniform with respect to $R$.
\par

\subsection{Technical Preliminaries }
In this section we outline certain energy and multiplier estimates, as well as estimates on the von Karman nonlinearity, which will be necessary in the proof of Theorem \ref{th:main2} below. For full details of the proof of Theorem \ref{th:main2} see \cite{delay}.
\smallskip\par
The following theorem is a case specialization found in \cite[Section 1.4, pp.38-45; Section 9.4, pp.496-497]{springer}; it is of critical importance in the analysis of von Karman plates without rotational inertia. The first  bound elucidates the local Lipschitz (quasi-Lipschitz) character of the von Karman nonlinearity, and is relatively recent and critical to our nonlinear analysis. The second bound is related to a compensated compactness of the nonlinear term, which by itself is also at the critical level. This  compensated compactness estimate is essential in proving ``smoothness" of  the attracting set.
\begin{theorem}\label{nonest}
Let $u^i \in \mathscr{B}_R(H^2_0(\Omega))$, $i=1,2$, and $z=u^1-u^2$.
 Then for $f(u) = - [u,v(u)+F_0]$ we have
\begin{equation}\label{f-est-lip}
||f(u^1)-f(u^2)||_{-\delta}  \le C_{\delta}\big(1+||u^1||_2^2+||u^2||_2^2\big)||z||_{2-\delta} \le C(\delta,R)||z||_{2-\delta}~~~ for~ all ~~\delta \in [0,1].
\end{equation}
If we further assume that $u^i \in C(s,t;H^2(\Omega))\cap C^1(s,t;L_2(\Omega))$, then we have  that
\begin{equation*}
- <f(u^1)-f(u^2),z_{t}> =\dfrac{1}{4}\frac{d}{dt}Q(z)+\frac{1}{2}
P(z)
\end{equation*}
where
\begin{equation*}
Q(z)=<v(u^1)+v(u^2),[z,z]> -||\Delta v(u^1+u^2,z)||^2
\end{equation*}
and
\begin{equation}\label{4.9aa}
P(z)=-<u^1_{t},[u^1,v(z)]> -<u^2_{t},[u^2,v(z)]> -<
u^1_{t}+u^2_{t},[z,v(u^1+u^2,z)]> .
\end{equation}
Moreover \cite{springer},
\begin{align}\label{eq4.5}
\Big|\int_s^t <f(u^1(\tau))-f(u^2(\tau)),z_t(\tau)>d\tau\Big| \le &~C(R)\sup_{\tau \in [s,t]} ||z||^2_{2-\eta}+\frac12 \Big|\int_s^tP(z)d\tau\Big|
\end{align} for some $0<\eta<1/2$ provided
 $u^i(\tau) \in \mathscr{B}_R(H^2_0(\Omega))$ for all $\tau\in [s,t].$
\end{theorem}

The above estimates and the standard displacement multiplier allow us to obtain 
the following estimates:
\begin{lemma}
Let $u^i \in C(0,T;\cH)\cap C^1(0,T;L_2(\Omega)) \cap L_2(-t^*,T;\cH)$ solve (\ref{plate}) with clamped or clamped-hinged boundary conditions and appropriate initial conditions on $[0,T]$ for $i=1,2$. Then the following estimate holds for all $\e >0$, for some $\eta > 0$, and $0 \le t \le T$:
\begin{align*}
\int_0^t \big(||\Delta u||^2-||u_t||^2\big)d\tau \le~\epsilon \int_0^t||u||^2_2 d\tau + C\int_{-t^*}^0 ||u(\tau)||_2^2 d\tau +\\ C(\epsilon,t^*,T)\sup_{\tau \in [0,t]} ||u(\tau)||^2_{2-\eta} + \e\\ -\int_0^t<f(u),u> d\tau + |<u_t(t), u(t)>| +|<u_t(s), u(s)>|.
\nonumber
\end{align*}
Moreover, in the case where we are considering the difference $z=u^1-u^2$ of solutions solving (\ref{difference}) with $u^i(t) \in \mathscr B_R(H^2(\Omega))$ for all $t\in [0,T]$, we may utilize the estimates in Theorem~\ref{nonest} (which eliminates the stand-alone $\e$) arrive at
\begin{align}\label{zmult}
\int_s^t \big(||\Delta z||^2-||z_t||^2\big)d\tau
\le &~\epsilon \int_s^t ||z||_2^2 d\tau +C\int_{s-t^*}^t ||z(\tau)||_{2-\s}^2 d\tau + \notag \\
C(\epsilon,T,R)\sup_{\tau \in [0,t]} ||z(\tau)||^2_{2-\eta} 
&+ E_z(t)+E_z(s),
\end{align}
where 
$E_z(t) \equiv \dfrac{1}{2}\big\{||\Delta z(t)||^2 + ||z_t(t)||^2\big\}$.
\end{lemma}
The final class of estimates we  need are energy estimates for the $z$ term defines as the solution 
to (\ref{difference}). Energy estimates for single solutions (making use of the nonlinear potential energy) 
can be derived straightforwardly from (\ref{energyrelation2}). 
The energy estimate on $z$, along with the estimate in (\ref{zmult}) above, will be used in 
showing asymptotic smoothness for the system.

Using the energy relation for the difference of two trajectories, \eqref{zmult}, 
and the integration by parts formula for the integral with delayed
term
we arrive at the following assertion.

\begin{lemma}\label{le:observbl}
Let $u^i \in C(0,T;\cH)\cap C^1(0,T;L_2(\Omega)) \cap L_2(-t^*,T;\cH)$ solve (\ref{plate}) with clamped boundary conditions and appropriate initial conditions on $[0,T]$ for $i=1,2$, $T\ge 2t^*$. Additionally, assume $u^i(t) \in \mathscr B_R(H^2(\Omega))$ for all $t\in [0,T]$.
Then the following estimates
\begin{align}\label{enest1}
 \frac{T}2 \Big[\Ez(T)+ &\int_{T-t^*}^T \Ez(\tau) d\tau\Big] \le ~ a_0\left(\Ez(0)+\int_{-t^*}^0 ||z(\tau)||_{2}^2 d\tau\right)
  +C(T,R)\sup_{\tau \in [0,T]}||z||^2_{2-\eta_*}
\\\nonumber &-a_1\int_0^Tds \int_s^T <f(u^1)-f(u^2),z_t>d\tau  -a_2\int_0^T <f(u^1)-f(u^2),z_t>d\tau
\end{align}
hold with $a_i$ independent of $T$ and $R$.
\end{lemma}

\subsection{Asymptotic Smoothness}
Recall that our dynamical system is $(S_t,\bH)$, where $S_t$ is the evolution operator corresponding 
to plate solutions to (\ref{plate}) and $\bH = H_0^2(\Omega) \times L_2(\Omega) \times L_2(-t^*,0;H_0^2(\Omega))$. 
To show asymptotic smoothness of this dynamical system, we will make use of some  abstract theorem, 
which generalizes the result given in \cite{kh}. In fact, criticality of nonlinear term prevents the  
use of other more standard methods, such as those given in  \cite{Babin-Vishik,lad,temam}.

To make use of this theorem, we will consider our functional $\Psi$ comprising lower order terms 
(compact with respect to $\bH$) and quasicompact ($\int_s^t <f(u^1)-f(u^2),z_t>d\tau$) terms. 
We need to produce an estimate which bounds trajectories in 
$\bH$, i.e. $||(u(t),u_t(t),u^t)||^2_{\bH}$ (taking the metric $d$ to be $||\cdot||_{\bH}$). 
Such an estimate will be produced below by combining our energy estimates produced earlier:
\begin{lemma}\label{le:khan}
Suppose $z=u^1-u^2$ is as in (\ref{difference}), with $y^i(t)=(u^i(t),u_t(t)^i,u^{t,i})$ and 
$y^i(t) \in \mathscr B_R(\bH)$ for all $t\ge 0$. Also, let $\eta >0$.
Then for every $0<\e<1$ there exists  $T=T_\e(R)$ such that the following estimate
\begin{equation*}
E_z(T) + \int_{T-t^*}^T ||z(\tau)||_2^2 d\tau \le \epsilon + \Psi_{\epsilon,T,R}(y^1,y^2)
\end{equation*}
 holds,
where \begin{align*} \Psi_{\epsilon,T,R}(y^1,y^2) \equiv& C(R,T) \sup_{\tau \in [0,T]} ||z(\tau)||_{2-\eta}^2 +a_1\Big| \int_0^T<f(u^1(\tau))-f(u^2(\tau)),z_t(\tau)>d\tau\Big| \\\nonumber& + a_2\Big|\int_0^T \int_s^T <f(u^1(\tau))-f(u^2(\tau)),z_t(\tau)>d\tau ds\Big|.\end{align*}
\end{lemma}
\begin{proof} It follows from \eqref{enest1} by dividing by $T$ and
taking $T$ large enough.
\end{proof}
In Lemma \ref{le:khan} above, we have obtained the necessary estimate for asymptotic smoothness; 
it now suffices to show that $\Psi$, as defined above, has the compensated compactness condition.

\begin{theorem}\label{smoothness}
The dynamical system $(S_t,\bH)$ generated by weak solutions to (\ref{plate}) is asymptotically smooth.
\end{theorem}
\begin{proof}
In line with the discussion above, we make use of 
the useful criterion (inspired by \cite{kh} and proven in \cite{ch-l}, see also
\cite[Chapter 7]{springer}) which reduces
asymptotic smoothness to finding a suitable functional on the state space
with a compensated compactness condition.
To do so, it suffices to show the compensated compactness condition for $\Psi_{\epsilon,T,R}$ which we now write as $\Psi$, with $\epsilon, T,$ and $R$ fixed along with the other constants given by the equation. Let $B$ be a bounded, positively invariant set in $\bH$, and let $\{y_n\}\subset B \subset \mathscr{B}_R(\bH)$. We must show that
$$
\liminf_m \liminf_n \Psi(y_n,y_m) = 0.
$$
The proof of this condition for $\Psi$ follows from the decomposition \cite[pp. 598--599]{springer}. 

\begin{align*}
<\mathcal{F}(z)(\tau) ,z_t(\tau)>=& ~\frac{1}{4} \frac{d}{d\tau} \big\{ - ||\Delta v(u^1) ||^2
- ||\Delta v(u^2) ||^2 + 2 <[z,z],F_0> \big\}\\\nonumber&-< [ v(u^2),u^2], u^1_t> - 
< [ v(u^1) , u^1 ], u^2_t>.
\end{align*}
Integrating the above expression in time and utilizing the sharp regularity of the Airy stress function, and the compactness of the lower order terms, the desired limit property follows.
For details see \cite{delay}.
\end{proof}
Having shown the asymptotic smoothness property, we can now conclude by Theorem \ref{dissmooth}
and Corollary \ref{c:a}  that there exists a compact global attractor $\mathbf{A} \subset \bH$ 
for the dynamical system $(S_t, \bH)$.
One may also see that this attractor has finite fractal dimension (see the discussion below).

As the final step in the proof of Theorem \ref{th:main2} we note that we have rewritten the dynamical system  generated by  the full flow-plate system (\ref{flowplate}) as the delayed  system in $(S_t, \bH)$ (\ref{reducedplate}). This is possible for sufficiently large times by Theorem \ref{rewrite}. We apply this result  to the dynamical system generated by the weak solution to (\ref{reducedplate}) on the space $\bH=\cH(\Omega)\times L_2(\Omega) \times L_2(-t^*,0;H_0^2(\Omega))$. This yields a compact global attractor $\mathbf{A} \subset \bH$ of finite dimension and additional regularity; to arrive at our set $\mathscr U$ (as in Theorem \ref{th:main2}) we then take $\mathscr U$ to be the projection of $\mathcal{A}$ on $\mathcal H \times L_2(\Omega)$, which concludes the proof.

\begin{remark}
{\rm
It should also be noted here that because we have rewritten our problem \eqref{flowplate} as a reduced delayed plate, and additionally changed the state space upon which we are operating, the  results obtained  on long-time behavior \textit{ will not be invariant } with respect to the flow component of the model, i.e. our global attractors will be with respect to the state space $\bH$, as defined above. Again, the data in the form of the delayed  term $u|_{(-t^*,0)}$
contains the information  from the flow itself. Obtaining global attractors for the full state space corresponding to $(u,u_t;\phi,\phi_t) \in H^2_0(\Omega) \times L_2(\Omega) \times H^1(\realsthree_+)\times L_2(\realsthree_+)$ is not  a
realistic task from  the  mathematical point of view. There is no damping imposed  on the system, thus the flow component
evolves according to  the full half space, unconstrained dynamics. The obtained result on the structure (without damping)  seems to be the best  possible result with respect to both the underlying physics and mathematics of the problem.
}
\end{remark}

\subsection{Quasistability Estimate and the Completion of the Proof of Theorem \ref{th:main2}}
In this section we refine our methods in the asymptotic smoothness calculation and work on 
trajectories from the attractor, whose existence has been established in the previous sections.

The key role is played by a {\em quasistability} estimate which reflects the fact that the flow can be stabilized exponentially
to a compact set. Alternatively, we might say that the flow is exponentially
stable, modulo a compact perturbation (lower order terms). The quadratic nature of 
the lower order terms is important for our considerations.

\begin{lemma}\label{quasi2}
Suppose that $u^1(t), u^2(t) \in \mathbf{A}$ are two solutions to \eqref{plate} with $z=u^1-u^2$. Then for $y=(z(t),z_t(t),z^t)$ there exists constants $\sigma, C(\sigma,\mathbf A)$ and $C$ such that for all $t>0$:
\begin{equation}\label{stab-estm}
||y(t)||_{\bH}^2 \le C(\sigma,\mathbf A)||y(0)||_{\bH}^2e^{-\sigma t}+C \sup_{\tau \in [0,t]} ||z(\tau)||_{2-\eta}^2.
\end{equation}
\end{lemma}
The estimate above is often referred to (in practice) as the ``quasistability" or ``stabilizability" 
estimate (see, e.g., \cite{chlJDE04,ch-l,springer}). 
\begin{proof}
Analyzing \eqref{enest1}, we may also write
\begin{align}\label{eq-obsr1}
T\left[E_z(T)+\int_{T-t^*}^TE_z(\tau)d\tau\right] \le c\big(E_z(0)+\int_{-t^*}^0||z(\tau)||_2^2d\tau\big)
\\
+
C\cdot T \sup_{s\in [0,T]}\Big|\int_s^T<\cF(z),z_t>d\tau \Big|+C(R,T)\sup_{\tau \in [s,t]} ||z||^2_{2-\eta},
\notag
\end{align}
where $\cF(z)= f(u^1)-f(u^2)$.
We note that $c$ does not depend on $T\ge \min\{1,2t^*\}$, and $l.o.t.$.

In order to prove the quasistability estimate, we have to handle the non-compact term $<
\mathcal{F}(z),z_{t}> $. We recall the  relation \eqref{eq4.5}  in Theorem \ref{nonest}:
if  $u^i \in C(s,t;H^2(\Omega))\cap C^1(s,t;L_2(\Omega))$ with $u^i(\tau) \in \mathscr{B}_R(L_2(\Omega))$
for $\tau\in [s,t]$, then
\begin{align}\label{1new}
\Big|
\int_s^t <\cF(z),z_t(\tau)>d\tau \le &~C(R)\sup_{\tau \in [s,t]} ||z||^2_{2-\eta}+\frac12\Big|
 \int_s^t P(z(\tau))d\tau\Big|
 \end{align} for some $0<\eta<1/2$.
Here $P(z)$ is given by (\ref{4.9aa}).
\par
Let $\gamma_{u^1} =\{(u^1(t),u^1_{t}(t),[u^1]^t):t\in \mathbb{R\}}$ and $\gamma_{u^2} =\{(u^2(t),u^2_{t}(t),[u^2]^t):t\in \mathbb{R\}}$ be trajectories from the
attractor $\mathbf{A}$. It is clear that for the pair $
u^1(t)$ and $u^2(t)$ satisfy the hypotheses of the estimate in \eqref{1new} for every
interval $[s,t]$.
Our main goal is to handle the second term on the right hand side of (\ref{1new}) which is of {\it critical
regularity}.
 To accomplish this we  use the already established {\it compactness} of the attractor in 
 the state space $\bH=\cH\times L_2(\Omega) \times L_2(-t^*,0;\cH)$ and also the method presented
 in \cite{springer}. This allows us
finally to conclude that for $y=(z(t),z_t(t),z^t)$ the estimate in \eqref{stab-estm}
is satisfied.
\end{proof}

The quasistability established in Lemma \ref{quasi2} makes it possible to
conclude that $\mathbf{A}$ has a finite fractal dimension.
\par
Additionally (see \cite{ch-l} or \cite{springer} for details) 
using quasistability estimate on a trajectory from $\mathbf{A}$ and a time shifted copy of the trajectory 
we can conclude that 
 $$ 
||u_{tt}(t)||^2+||u_t(t)||_2^2 +||u(t)||_{4}^{2}\leq C ~\text{for all} ~
t\in \reals.$$
on the attractor.
Thus, we have additional regularity of the trajectories from the attractor $\mathbf A \subset \bH$
stated in Theorem~\ref{th:main2}.
\par

\end{document}